\documentclass[letterpaper,11pt]{article}
\usepackage{latexsym,amssymb,enumerate,amsmath,epsfig,amsthm,dsfont,bm}
\usepackage[margin=1in]{geometry}
\usepackage{setspace,color}
\usepackage{tikz}
\usepackage{floatrow}
\usepackage{graphicx,subfigure}
\usepackage[ruled]{algorithm2e}
\usepackage{algorithmic}
\usepackage{epstopdf}
\usepackage{grffile}

\usepackage{hyperref}
\hypersetup{
	colorlinks=true,
	linkcolor=blue,
	citecolor=blue,
	filecolor=blue,      
	urlcolor=blue
}
\usepackage{lineno}

\newcommand{\zero}{\mathbf{0}}
\newcommand{\x}{\mathbf{x}}

\newcommand{\p}{\mathbf{p}}
\newcommand{\A}{\mathbf{A}}
\newcommand{\M}{\mathbf{M}}
\newcommand{\bP}{\mathbf{P}}
\newcommand{\bS}{\mathbf{S}}

\newcommand{\cK}{\mathcal{K}}

\newcommand{\n}{\mathbf{n}}
\newcommand{\D}{\mathbf{D}}
\newcommand{\I}{\mathbf{I}}
\newcommand{\RR}{\mathds{R}}
\newcommand{\cT}{\mathcal{T}}
\newcommand{\cof}{\mathrm{cof}}
\newtheorem{remark}{Remark}[section]
\newtheorem{theorem}{Theorem}[section]

\DeclareMathOperator*{\argmin}{arg\,min}
\begin{document}

\title{An Efficient Operator--Splitting Method for the Eigenvalue Problem of the Monge-Amp\`{e}re Equation}
\date{\large{\emph{In memory of Professor Roland Glowinski}}}
\author{ Hao Liu \thanks{Department of Mathematics, Hong Kong Baptist University, Kowloon Tong, Hong Kong (Email: {\bf haoliu@hkbu.edu.hk}).}
	\and
	Shingyu Leung\thanks{Department of Mathematics, the Hong Kong University of Science and Technology, Clear Water Bay, Hong Kong (Email: {\bf masyleung@ust.hk})}
 \and
  Jianliang Qian \thanks{Department of Mathematics and Department of CMSE, Michigan State University, East Lansing, MI 48824 (Email: {\bf jqian@msu.edu}).}
}

\maketitle
\begin{abstract}
We develop an efficient operator--splitting method for the eigenvalue problem of the Monge--Amp\`{e}re operator in the Aleksandrov sense. The backbone of our method relies on a convergent Rayleigh inverse iterative formulation proposed by Abedin and Kitagawa (Inverse iteration for the {M}onge--{A}mp{\`e}re eigenvalue problem, {\it Proceedings of the American Mathematical Society}, 148 (2020), no. 11, 4975-4886). Modifying the theoretical formulation, we develop an efficient algorithm for computing the eigenvalue and eigenfunction of the Monge--Amp\`{e}re operator by solving a constrained Monge--Amp\`{e}re equation during each iteration. Our method consists of four essential steps: (i) Formulate the Monge--Amp\`{e}re eigenvalue problem as an optimization problem with a constraint; (ii) Adopt an indicator function to treat the constraint; (iii) Introduce an auxiliary variable to decouple the original constrained optimization problem into simpler optimization subproblems and associate the resulting new optimization problem with an initial value problem; and (iv) Discretize the resulting initial-value problem by an operator--splitting method in time and a mixed finite element method in space. The performance of our method is demonstrated by several experiments. Compared to existing methods, the new method is more efficient in terms of computational cost and has a comparable rate of convergence in terms of accuracy.

\end{abstract}

\section{Introduction}
\label{sec.intro}
The Monge-Amp\`{e}re equation is a second-order fully nonlinear PDE in the form of 
\begin{align}
		\det \D^2u=f,
		\label{eq.MA}
\end{align}
where $\D^2u$ denotes the Hessian of $u$. The Monge-Amp\`{e}re equation originates from differential geometry in which it describes a surface with prescribed Gaussian curvature \cite{bakelman2012convex,kazdan1985prescribing}. The existence, uniqueness and regularity of the solution has been  extensively studied \cite{bakelman2012convex,roberts1995fully,gilbarg1977elliptic}, and related applications can be found in optimal transport \cite{benamou2000computational,froese2012numerical}, seismology \cite{engquist2016optimal}, image processing \cite{haker2004optimal}, finance \cite{stojanovic2004optimal}, and  geostrophic flows \cite{feng2009modified}. 

Due to its broad applications, in the past decade, a lot of efforts have been devoted to developing numerical methods for the Monge-Amp\`{e}re equation. One line of research is to develop wide-stencil based finite-difference schemes \cite{froese2011convergent,froese2011fast} for equation (\ref{eq.MA}) with Dirichlet boundary conditions. Such a class of methods utilizes the fact that $\det \D^2u$ equals the product of the eigenvalues of $\D^2u$, so that these methods use wide-stencils to estimate the eigenvalues. Later on, such methods were extended to accommodate transport boundary conditions in \cite{froese2012numerical}. Another line of research is to design finite-element based methods. In \cite{feng2009vanishing, feng2009mixed}, the authors  proposed the vanishing moment method, which approximates a fully nonlinear second-order PDE by a fourth-order PDE. In \cite{dean2003numerical, dean2004numerical,caboussat2013least,caboussat2018least}, the authors formulate equation (\ref{eq.MA}) as an optimization problem. Fast augmented Lagrangian algorithms are then designed to solve the new problems. Recently, operator--splitting methods have been proposed in \cite{glowinski2019finite, liu2019finite}. Taking advantage of the divergence form of $\det\D^2u$, the authors of \cite{glowinski2019finite, liu2019finite} decouple the nonlinearity of equation (\ref{eq.MA}) by introducing an auxiliary variable so that solving equation (\ref{eq.MA}) is reduced to finding the steady-state solution of an initial value problem, which is time-discretized by an operator--splitting method and space-discretized by a mixed finite-element method. Other numerical methods for equation (\ref{eq.MA}) include \cite{awanou2015standard, benamou2010two,benamou2014numerical,feng2022narrow,caboussat2021second,caboussat2021adaptive}; see the survey \cite{feng2013recent} for more related works. 

Existing works discussed above target equation (\ref{eq.MA}) with various boundary conditions. Another interesting problem of the Monge--Amp\`{e}re type is the eigenvalue problem, reading as
\begin{align}
	\begin{cases}
		\det(\D^2u)=\lambda |u|^d &\mbox{ in } \Omega,\\
		u=0 &\mbox{ on } \partial \Omega,
	\end{cases}
	\label{eq.MAEV}
\end{align}
where $\Omega\subset \RR^d$ ($d\geq 2$) is an open bounded convex domain, and $\lambda=\lambda[\Omega]$ is the unknown eigenvalue of the Monge--Amp\`{e}re operator on $\Omega$. Problem (\ref{eq.MAEV}) was first studied by Lions in \cite{lions1985two} and later by Tso in \cite{tso1990real}. They proved the existence, uniqueness and regularity of the solution on an open, bounded, smooth, uniformly convex domain. The result was then extended by Le in \cite{le2017eigenvalue} to general  bounded convex domains. Theoretically, to find the solution of equation (\ref{eq.MAEV}), a variational formulation was proposed in \cite{tso1990real}, and a convergent Rayleigh quotient inverse iterative formulation was proposed in \cite{abedin2020inverse} which was further improved in \cite{le2020convergence}. Since, during each Rayleigh quotient iteration, the algorithm in \cite{abedin2020inverse} requires solving a Monge--Amp\`{e}re type equation, how to efficiently implement this formulation numerically has not been studied. The only work on the numerical solution of equation (\ref{eq.MAEV}) we are aware of is \cite{glowinski2020numerical}, in which the authors proposed operator--splitting methods for a class of Monge-Amp\`{e}re eigenvalue problem. In \cite{glowinski2020numerical}, taking advantage of the divergence form, the authors takes equation (\ref{eq.MAEV}) as the optimality condition of a constrained optimization problem, in which $\lambda$ is considered as the Lagrange multiplier, and an operator--splitting method was proposed to solve the new problem. 

Similar to equation (\ref{eq.MA}), the eigenvalue problem (\ref{eq.MAEV}) is a fully nonlinear second-order PDE. One effective way to solve such PDEs is the operator--splitting method, which decomposes complicated problems into several easy--to--solve subproblems by introducing auxiliary variables. Then the new problem will be formulated as solving an initial value problem, which is then time discretized using operator--splittings. All variables will be updated in an alternative fashion, where each subproblem either has an explicit solution or can be solved efficiently. The operator--splitting method has been applied to numerically solving PDEs \cite{glowinski2019finite,liu2019finite}, image processing \cite{liu2021color,deng2019new,liu2021operator,liu2022elastica}, surface reconstruction \cite{he2020curvature}, inverse problems \cite{glowinski2015penalization}, obstacle problems \cite{liu2022fast}, and computational fluid dynamics \cite{bukavc2013fluid,bonito2016operator}. We refer readers to monographs \cite{glowinski2017splitting,glowinski2016some} for detailed discussions on operator--splitting methods.

In this work, we propose an efficient numerical implementation of the formulation proposed in \cite{abedin2020inverse} to compute the eigenvalue and eigenfunction of the Monge--Amp\`{e}re operator on an open, bounded, convex domain $\Omega$. Since each Rayleigh quotient inverse iteration of the formulation in \cite{abedin2020inverse} requires solving a Monge--Amp\`{e}re equation, we first use the divergence form of the Monge--Amp\`{e}re operator to rewrite the problem as an optimization problem. To stabilize our formulation, we consider a constrained version of the optimization problem by forcing the eigenfunction $u$ to have unit $L_2$-norm: $\|u\|_2=1$. The constrained problem is converted to an unconstrained problem by utilizing an indicator function of the constraint set. Then we decouple the nonlinearity of the functional by introducing an auxiliary variable, and we associate it with an initial value problem in the flavor of gradient flow. The initial value problem is time discretized by an operator-splitting method and space discretized by a mixed finite-element method in the space of piecewise-linear continuous functions. The efficiency of the proposed method is demonstrated by several numerical experiments.

We organize the rest of this article as follows: We introduce the background and  summarize the convergent formulation of \cite{abedin2020inverse} for equation (\ref{eq.MAEV}) in Section \ref{sec.background}. Our new operator-splitting approach for implementing this convergent formulation is presented in Section \ref{sec.formulation}. Our operator-splitting scheme is time discretized in Section \ref{sec.os} and space discretized in Section \ref{sec.discretization}. We demonstrate the efficiency of the proposed method by several numerical experiments in Section \ref{sec.experiments} and conclude this article in Section \ref{sec.conclusion}.


\section{A convergent inverse iteration for the eigenvalue problem}
\label{sec.background}

Let $\Omega\subset \RR^d$ be an open bounded convex domain. In equation (\ref{eq.MAEV}), if $u$ is a convex function, one has $u\leq 0$ and $|u|=-u$. The existence and uniqueness of the eigen-pair was studied in \cite{lions1985two}:
\begin{theorem}
	Assume that $\Omega\subset \RR^d$ is a smooth, bounded, uniformly convex domain. There exist a unique positive constant $\lambda_{\rm MA}$ and a unique (up to positive multiplicative constants) nonzero convex function $u\in C^{1,1}(\bar{\Omega})\cap C^{\infty}(\Omega)$ solving the eigenvalue problem (\ref{eq.MAEV}). The constant $\lambda_{\rm MA}$ is called the Monge-Amp\`{e}re eigenvalue of $\Omega$ and $u$
is called a Monge-Amp\`{e}re eigenfunction of $\Omega$. 	
\end{theorem}

Define the Rayleigh quotient of a function $u$ for the Monge-Amp\`{e}re operator as
\begin{align}
	R(u)=\frac{\int_{\Omega} -u\det(\D^2u)d\x}{\int_{\Omega}(-u)^{d+1}d\x},
	\label{eq.Rayleigh}
\end{align}
and the function space $\cK$ as
\begin{align*}
	\cK=\left\{u\in C^{0,1}(\bar{\Omega})\cap C^{\infty}(\Omega): u \mbox{ is convex and nonzero in } \Omega, \ u=0 \mbox{ on } \partial\Omega \right\}.
\end{align*}

Tso \cite{tso1990real} showed that $\lambda_{\rm MA}$ can be written as the infimum of Rayleigh quotients:
\begin{theorem}
	Assume that $\Omega\subset \RR^d$ is a smooth, bounded and uniformly convex domain. Then
	\begin{align}
		\lambda_{\rm MA}=\inf_{u\in \cK} R(u).
		\label{eq.Ray}
	\end{align}
\label{thm.eigen}
\end{theorem}

Based on the property (\ref{eq.Ray}), the following inverse iterative scheme for the eigenvalue problem (\ref{eq.MAEV}) was proposed by Abedin and Kitagawa in \cite{abedin2020inverse}:
\begin{align}
	\begin{cases}
		u^0=u_0,\\
		\det(\D^2 u^{k+1})=R(u^k)|u^k|^d & \mbox{ in } \Omega,\\
		u^{k+1}=0 &\mbox{ on } \partial \Omega,
	\end{cases}
\label{eq.scheme}
\end{align}
where $u_0$ is a given initial condition, and they further proved the convergence of the inverse iteration:
\begin{theorem}
	Assume that $\Omega\subset \RR^d$ is an open bounded convex domain. Let $u_0\in C(\bar{\Omega})$ satisfy the following: 
	\begin{enumerate}
		\item[(i)] $u_0$ is convex and $u_0\leq 0$ on $\partial \Omega$;
		\item[(ii)] $R(u_0)<\infty$; 
		\item[(iii)] $\det(\D^2 u_0) \geq c_0$ in $\Omega$, where $c_0$ is some positive constant.
	\end{enumerate}
Then, for $k>0$, $u^k$ in equation (\ref{eq.scheme}) converges uniformly on $\bar{\Omega}$ to a nonzero Monge-Amp\`{e}re eigenfunction, and $R(u_k)$ converges to $\lambda_{\rm MA}$.
\label{thm.convergence}
\end{theorem}

Theorem \ref{thm.convergence} was improved in \cite{le2020convergence} so that conditions (i) and (iii) are removed; consequently, the inverse iteration  converges for all convex initial data having finite and nonzero Rayleigh quotient to a nonzero Monge-Amp\`{e}re eigenfunction of $\Omega$.

\section{A modified formulation of the inverse iteration}
\label{sec.formulation}

Given an initial convex function $u_0$ with bounded nonzero Rayleigh quotient, the inverse iteration (\ref{eq.scheme}) generates the sequence 
$\{(R(u^k),u^k)\}$ which is guaranteed to converge to the solution of the eigenvalue problem (\ref{eq.MAEV}). When updating $u^{k+1}$ from $u^k$, one needs to solve a Monge-Amp\`{e}re equation with the Dirichlet boundary condition, which is a nonlinear problem. It has not been studied yet how to implement the inverse iteration efficiently to produce numerical approximations to the eigenvalue problem of the Monge-Amp\`{e}re operator.
Therefore, we are motivated to develop an efficient algorithm to implement this inverse iterative method. 

To achieve this purpose, we adopt a recently developed operator-splitting method (see \cite{glowinski2019finite, liu2019finite, glowinski2020numerical}) to solve equation (\ref{eq.scheme}) numerically. We focus on the case $d=2$. Our method can be easily extended to higher dimensional problems.

We first reformulate equation (\ref{eq.scheme}) using the following identity:
\begin{align}
	\det(\D^2u)=\frac{1}{2}\nabla\cdot (\cof(\D^2u) \nabla u),
	\label{eq.div}
\end{align}
where 
$
\cof(\D^2u)=
\begin{bmatrix}
	\frac{\partial^2 u}{\partial x_1^2} & -\frac{\partial^2 u}{\partial x_1 \partial x_2}\\
	-\frac{\partial^2 u}{\partial x_1 \partial x_2} & \frac{\partial^2 u}{\partial x_2^2}
\end{bmatrix}
$
is the cofactor matrix of $\D^2u$.

Incorporating equation (\ref{eq.div}) into equations (\ref{eq.scheme}) and (\ref{eq.Rayleigh}) gives rise to
\begin{align}
	\begin{cases}
		u^0=u_0,\\
		\nabla\cdot (\cof(\D^2u^{k+1}) \nabla u^{k+1})={2}\;R(u^k)|u^k|^d & \mbox{ in } \Omega,\\
		u^{k+1}=0 &\mbox{ on } \partial\Omega,
	\end{cases}
	\label{eq.scheme.div}
\end{align}
with 
\begin{align}
	R(u)=\frac{\displaystyle\int_{\Omega} (\cof(\D^2u) \nabla u)\cdot \nabla u d\x}{2\displaystyle\int_{\Omega}(-u)^{3}d\x},
	\label{eq.Rayleigh.div}
\end{align}
where we used integration by parts when deriving equation (\ref{eq.Rayleigh.div}). 

From equation \eqref{eq.scheme.div}, updating $u^{k+1}$ from $u^k$ is equivalent to solving the optimization problem
\begin{align}
	\begin{cases}
		\min\limits_{w}\displaystyle \left[\int_{\Omega} (\cof(\D^2w) \nabla w)\cdot \nabla wd\x+6\displaystyle\int_{\Omega} f^kwd\x\right],\\
		w=0 \mbox{ on } \partial\Omega,
	\end{cases}
	\label{eq.min0}
\end{align}
with $f=R(u^k)|u^k|^2$, which can be derived from the first-order variational principle; see \cite{glowinski2019finite, liu2019finite}. Note that if $(\lambda_{\rm MA},u^*)$ is a solution to equation (\ref{eq.MAEV}), $(\lambda_{\rm MA},\alpha u^*)$ is also a solution for any $\alpha> 0$ (assuming that we are looking for convex eigenfunctions). To make the solution of equation (\ref{eq.MAEV}) unique, we restrict our attention to looking for the eigenfunction $u^*$ satisfying 
\begin{align}
	\|u^*\|_2=1.
	\label{eq.constraint}
\end{align}
Therefore it is natural to add the constraint $\|w\|_2=1$ to equation (\ref{eq.min0}). However, usually a constrained optimization problem is more challenging to solve than an unconstrained one. Therefore, to remove the constraint while enforcing $\|w\|_2=1$, we utilize an indicator function.


Define the set
\begin{align*}
	S=\{w: w \mbox{ is smooth}, \|w\|_2=1\}
\end{align*}
and its indicator function
\begin{align*}
	I_{S}(w)=\begin{cases}
		0 &\mbox{ if } w\in S,\\
		+\infty & \mbox{ otherwise}.
	\end{cases}
\end{align*}

Equation (\ref{eq.min0}) with constraint $\|w\|_2=1$ can be rewritten as
\begin{align}
	\begin{cases}
		\min\limits_{w}\displaystyle \left[\int_{\Omega} (\cof(\D^2w) \nabla w)\cdot \nabla wd\x+6\displaystyle\int_{\Omega} f^k\,wd\x +I_{S}(w) \right],\\
		w=0 \mbox{ on } \partial\Omega.
	\end{cases}
	\label{eq.min}
\end{align}

We follow \cite{glowinski2019finite} to introduce a matrix-valued auxiliary variable $\p$ to decouple the nonlinearity in equation (\ref{eq.min}). Then solving equation (\ref{eq.min}) is equivalent to solving 
\begin{align}
	\begin{cases}
		\min\limits_{w, \p}&\displaystyle \left[\int_{\Omega} (\cof(\p) \nabla w)\cdot \nabla wd\x+6\displaystyle\int_{\Omega} f^k\,wd\x+I_{S}(w) \right],\\
		w=0 &\mbox{ on } \partial\Omega,\\
		\p=\D^2w & \mbox{ in } \Omega.
	\end{cases}
	\label{eq.min.1}
\end{align}

After computing the Euler-Lagrange equation, if $(v,\p)$ is a solution to equation (\ref{eq.min.1}), we have
\begin{align}
	\begin{cases}
		\nabla\cdot (\cof(\p) \nabla v)-2f^k+\partial I_{S}(v) \ni 0 & \mbox{ in } \Omega,\\
		v=0 &\mbox{ on } \partial\Omega,\\
		\p=\D^2v,&\mbox{ in } \Omega,	
         \end{cases}
	\label{eq.cMA.1}
\end{align}
where $\partial I_{S}$ denotes the sub-differential of $I_{S}$. 

We associate equation (\ref{eq.cMA.1}) with the following initial value problem (in the flavor of gradient flow)
\begin{align}
	\begin{cases}
		\begin{cases}
			\frac{\partial v}{\partial t}+\nabla\cdot (\varepsilon \I+ \cof(\p) \nabla v)-2f^k+\partial I_{S}(v)\ni 0 & \mbox{ in } \Omega\times(0,+\infty),\\
			v=0 &\mbox{ on } \partial\Omega\times (0,+\infty),
		\end{cases}\\
		\frac{\partial \p}{\partial t}+\gamma(\p-D^2v)=\mathbf{0} \quad \mbox{ in } \Omega\times(0,+\infty),\\
		v(0)=v_0,\ \p(0)=\p_0,
	\end{cases}
	\label{eq.cMA.2}
\end{align}
where $\I$ is the identity matrix, $\zero$ is the zero matrix, and $\varepsilon>0$ is a small constant. The term $\varepsilon \I$ is a regularization term in order to handle the case that $\inf_{\x\in\Omega} f^k(\x)=0$. Then $u^{k+1}$ is the steady state of $v$.

In equation (\ref{eq.cMA.2}), $\gamma$ controls the evolution speed of $\p$. A natural choice is to let $\p$ evolve with a similar speed as that of $v$, leading to 
\begin{align*}
	\gamma=\beta \lambda_0
\end{align*}
with $\lambda_0$ being the smallest eigenvalue of $-\nabla^2$ and $\beta>0$ being some constant. 

\section{An operator splitting method to solve equation (\ref{eq.cMA.2})}
\label{sec.os}
\subsection{The operator splitting strategy}
The structure of equation (\ref{eq.cMA.2}) is well--suited to be time-discretized by the operator splitting method. Among many possible discretization schemes, we choose the simplest Lie scheme. 

Let $\tau>0$ denote the time step and denote $t^n=n\tau$. We time-discretize equation (\ref{eq.cMA.2}) as follows:\\
Initialization:
\begin{align}
	v^0=v_0,\ \p^0=\p_0.
	\label{eq.split.0}
\end{align}
For $n>0$, update $(v^n,\p^n)\rightarrow(v^{n+1/3},\p^{n+1/3}) \rightarrow(v^{n+2/3},\p^{n+2/3})\rightarrow (v^{n+1},\p^{n+1})$ as:\\
\textbf{Step 1}: Solve
\begin{align}
	\begin{cases}
		\begin{cases}
			\frac{\partial v}{\partial t}+\nabla\cdot (\varepsilon \I+ \cof(\p) \nabla v)-2f^k= 0 & \mbox{ in } \Omega\times(t^n,t^{n+1}),\\
			v=0 &\mbox{ on } \partial\Omega\times (t^n,t^{n+1}),
		\end{cases}\\
		\frac{\partial \p}{\partial t}=\mathbf{0} \quad \mbox{ in } \Omega\times(t^n,t^{n+1}),\\
		v(t^n)=v^n,\ \p(t^n)=\p^n,
	\end{cases}
\label{eq.split.1}
\end{align}
and set $v^{n+1/3}=v(t^{n+1}), \ \p^{n+1/3}=\p(t^{n+1})$.\\
\textbf{Step 2}: Solve
\begin{align}
	\begin{cases}
		\begin{cases}
			\frac{\partial v}{\partial t}= 0 & \mbox{ in } \Omega\times(t^n,t^{n+1}),\\
			v=0 &\mbox{ on } \partial\Omega\times (t^n,t^{n+1}),
		\end{cases}\\
		\frac{\partial \p}{\partial t}+\gamma(\p-D^2v)=\mathbf{0} \quad \mbox{ in } \Omega\times(t^n,t^{n+1}),\\
		v(t^n)=v^{n+1/3},\ \p(t^n)=\p^{n+1/3},
	\end{cases}
\label{eq.split.2}
\end{align}
and set $v^{n+2/3}=v(t^{n+1}), \ \p^{n+2/3}=\p(t^{n+1})$.\\
\textbf{Step 3}: Solve
\begin{align}
	\begin{cases}
		\begin{cases}
			\frac{\partial v}{\partial t}+\partial I_{S}(v)\ni 0 & \mbox{ in } \Omega\times(t^n,t^{n+1}),\\
			v=0 &\mbox{ on } \partial\Omega\times (t^n,t^{n+1}),
		\end{cases}\\
		\frac{\partial \p}{\partial t}=\mathbf{0} \quad \mbox{ in } \Omega\times(t^n,t^{n+1}),\\
		v(t^n)=v^{n+2/3},\ \p(t^n)=\p^{n+2/3},
	\end{cases}
	\label{eq.split.3}
\end{align}
and set $v^{n+1}=v(t^{n+1}), \ \p^{n+1}=\p(t^{n+1})$.

The scheme (\ref{eq.split.0})--(\ref{eq.split.3}) is only semi-constructive since one still needs to solve the subproblems in equations (\ref{eq.split.1})--(\ref{eq.split.3}). For equation (\ref{eq.split.2}), we have the explicit solution for $\p^{n+2/3}$:
\begin{align*}
	\p^{n+2/3}=e^{-\gamma \tau}\p^n+(1-e^{-\gamma\tau})\D^2v^{n+1/3}.
\end{align*}
Since the solution of equation (\ref{eq.MAEV}) is a convex function,  the Hessian $\D^2u$ is a semi-positive definite matrix. Note that $\p$ is an auxiliary variable estimating $\D^2v$,  we project it onto the space of semi-positive definite symmetric matrices once $\p^{n+2/3}$ is computed. We denote the projection operator by $P_+$; see more details in Section \ref{sec.projP}. 

For other subproblems, we adopt the one-step backward Euler scheme (the Markchuk-Yanenko type). Our updating formulas are summarized as follows:
\begin{align}
	&\begin{cases}
		\frac{v^{n+1/3}-v^n}{\tau}+\nabla\cdot (\varepsilon \I+ \cof(\p^n) \nabla v^{n+1/3})-2f^k= 0 &\mbox{ in } \Omega,\\
		v^{n+1/3}=0 &\mbox{ on } \partial \Omega,
	\end{cases}
	\label{eq.alg.1}\\
	&\quad \p^{n+1}=P_+\left(e^{-\gamma \tau}\p^n+(1-e^{-\gamma\tau})\D^2v^{n+1/3}\right), \label{eq.alg.2}\\
	&\begin{cases}
		\frac{v^{n+1}-v^{n+1/3}}{\tau}+\partial I_{S}(v^{n+1})\ni 0 &\mbox{ in } \Omega,\\
		v^{n+1}=0 &\mbox{ on } \partial \Omega.
	\end{cases}
	\label{eq.alg.3}
\end{align}

\begin{remark}
	Equation (\ref{eq.cMA.2}) is very similar to problem (36) in \cite{glowinski2020numerical}, except that in our current scheme the constraint is $\|u\|_2=1$ and that in \cite{glowinski2020numerical} it is $\|u\|_{3}=1$. Despite similar formulations, the numerical treatments are very different. In equations (\ref{eq.alg.1})-(\ref{eq.alg.3}), $f^k$ and the indicator function $\partial I_{S}$ are separately distributed into two sub-steps. Equation (\ref{eq.alg.3}) simply results in a projection to the unit sphere; see Section \ref{sec.projection} for details.
	
	In \cite{glowinski2020numerical}, $\lambda\,d\,u|u|$ with $d$ being the spatial dimension plays the role of $f^k$ and the constraint plays the role of $\partial I_{S}$, and both terms are arranged in the same sub-step (problem (50b) in \cite{glowinski2020numerical}):
	\begin{align}
		\begin{cases}
			u^{n+2/3}-u^{n+1/3}=3\tau \lambda^{n+1} u^{n+2/3}|u^{n+2/3}|,\\
			{\displaystyle \int_{\Omega}} |u^{n+2/3}|^{3}d\x=1
		\end{cases}
	\label{eq.Lagrange}
	\end{align} 
The constraint $\|u\|_{3}=1$ cannot be replaced by $\|u\|_{2}=1$ since equation (\ref{eq.Lagrange}) was considered as an optimality condition of a Lagrangian functional and $\tau\lambda^{n+1}$ is the Lagrange multiplier. As a result, $u^{n+2/3}$ solves
\begin{align}
	u^{n+2/3}\in \argmin_{v:\int_{\Omega} |v|^{3}d\x=1} \left[ \frac{1}{2} \int_{\Omega} |v|^2d\x -\int_{\Omega} u^{n+1/3}vd\x\right].
	\label{eq.Lagrange.1}
\end{align}
Unlike (\ref{eq.alg.3}), the solution to problem (\ref{eq.Lagrange.1}) does not have an explicit expression, so that an iterative method (such as sequential quadratic programming) was used in \cite{glowinski2020numerical} to solve problem (\ref{eq.Lagrange.1}).
\end{remark}
\begin{remark}
	Compared to the algorithm (\ref{eq.scheme}) proposed in \cite{abedin2020inverse}, our scheme has an additional term related to the constraint $\|u\|_2=1$, and such a constraint leads to the projection step (\ref{eq.alg.3}) which helps stabilize our numerical algorithm.
\end{remark}
\subsection{On the solution to equation (\ref{eq.alg.3})}
\label{sec.projection}
In the scheme above, problems (\ref{eq.alg.1}) and (\ref{eq.alg.2}) are easy to solve. In equation (\ref{eq.alg.3}), $v^{n+1}$ solves
\begin{align}
	\begin{cases}
		\min\limits_w \left[\frac{1}{2\tau}\displaystyle\int_{\Omega} \|w-v^{n+1/3}\|_2^2 d\x +I_{S}(w)\right],\\
		w=0 \mbox{ on } \partial\Omega.
		\label{eq.v.min}
	\end{cases}
\end{align}
Since $I_{S}(w)$ is the indicator function of $S$ in which $\|w\|_2=1$, the exact solution of equation (\ref{eq.v.min}) reads as
\begin{align}
	v^{n+1}=\frac{v^{n+1/3}}{\|v^{n+1/3}\|_2}.
	\label{eq.v.proj}
\end{align}

\subsection{On the initial condition}
We next discuss the initial condition $u_0$ in the outer iteration and $(v_0,\p_0)$ in the inner iteration. The convergence theorem for the scheme (\ref{eq.scheme}), Theorem \ref{thm.convergence}, requires the initial condition to be convex and smooth. A simple choice is to set $u_0$ as the solution to
\begin{align}
	\begin{cases}
		\det \D^2 u_0=1 &\mbox{ in } \Omega,\\
		u_0=0 &\mbox{ on } \partial \Omega.
	\end{cases}
\label{eq.initial.0}
\end{align}
However, solving equation (\ref{eq.initial.0}) is not trivial. Since $u_0$ is only the initial condition and the iterates generated by the inverse iteration are eventually smooth as shown in \cite{le2020convergence}, we do not need to solve equation (\ref{eq.initial.0}) exactly. An operator splitting method is proposed in \cite{glowinski2019finite} to solve equation (\ref{eq.initial.0}). To make the initialization simpler, we will choose $u_0$ as the initial condition according to a strategy used in \cite{glowinski2019finite}. Specifically, $u_0$ is the solution to the Poisson problem
\begin{align}
	\begin{cases}
		\nabla^2 u_0=2\eta & \mbox{ in } \Omega,\\
		u_0=0 & \mbox{ on } \partial \Omega,
	\end{cases}
\label{eq.initial.outer}
\end{align}
where $\eta>0$ is of $O(1)$. 
 
For the initial condition $(v_0,\p_0)$ in the $k+1$-th outer iteration, we simply set 
\begin{align}
	v_0=u^k, \ \p_0=\D^2 v_0.
	\label{eq.initial.inner}
\end{align}

Our algorithm is summarized in Algorithm \ref{alg.1}.
\begin{algorithm}
	\caption{\label{alg.1}An operator-splitting method for solving problem (\ref{eq.MAEV})}
	\begin{algorithmic}
		\STATE {\bf Input:} Parameters $\gamma,\tau,\varepsilon, N$.
		\STATE {\bf Initialization:} Set $k=0$. Initialize $u^0$ according to equation (\ref{eq.initial.outer}). 
		\WHILE{not converge}
		\STATE Step 1. Compute $f^k=R(u^k)|u^k|^2$ according to equation (\ref{eq.Rayleigh.div}).
		\STATE Step 2. Set $n=0$. Initialize $(v^0,\p^0)$ according to equation (\ref{eq.initial.inner}).
		\WHILE{not converge}
		\STATE Step 3.1. Solve equation (\ref{eq.alg.1}) for $v^{n+1/3}$.
		\STATE Step 3.2. Solve equation (\ref{eq.alg.2}) for $\p^{n+1}$.
		\STATE Step 3.3. Solve equation (\ref{eq.alg.3}) for $v^{n+1}$.
		\STATE Step 3.4. Set $n=n+1$.
		\ENDWHILE
		\STATE Step 4. Set $u^{k+1}$ as the converged $v^*$.
		\STATE Step 5. Set $k=k+1$.
		\ENDWHILE
		\STATE {\bf Output:} The converged eigenfunction $u^*$ and eigenvalue $\lambda_{\rm MA}$.
	\end{algorithmic}
\end{algorithm}

\section{A finite element implementation of scheme (\ref{eq.alg.1})-(\ref{eq.alg.3})}
\label{sec.discretization}
\subsection{Generalities}
Let $\Omega\subset \RR^2$ be an open bounded convex polygonal domain (or it has been approximated by such a domain). Let $\cT_h$ be a triangulation of $\Omega$, where $h$ denotes the length of the longest edge of triangles in $\cT_h$. Define the following two piecewise linear function spaces
\begin{align*}
	&V_h=\{ \phi\in C^0(\bar{\Omega}): \phi_T\in \bP_1 \mbox{ for }\forall T\in \cT_h\},\\
	&V_{0h}=\{\phi\in V_h:\phi|_{\partial\Omega}=0\},
\end{align*}
where $\bP_1$ is the space of polynomials of two variables with degree no larger than 1. Let $H^1(\Omega)$ be the Sobolev space of order 1 and $H^1_0(\Omega)$ be the collection of functions in $H^1(\Omega)$ with vanishing trace on $\partial\Omega$. Then $V_h$ and $V_{0h}$ are approximations of $H^1(\Omega)$ and $H^1_0(\Omega)$, respectively.

Denote the set of vertices of $\cT_h$ by $\Sigma_h$. We further denote the interior vertices of $\cT_h$ by $\Sigma_{0h}=\Sigma_h\backslash(\Sigma_h\cap \partial\Omega)$. We use $N_h$ and $N_{0h}$ to denote the cardinality of $\Sigma_h$ and $\Sigma_{0h}$, respectively. We have
\begin{align*}
	\dim V_h=N_h \quad \mbox{ and } \quad \dim V_{0h}=N_{0h}.
\end{align*}

We order the vertices of $\cT_h$ so that $\Sigma_{0h}=\{Q_l\}_{l=1}^{N_{0h}}$, where $Q_l$'s denote the vertices. For any $1\leq l\leq N_h$, we use $\omega_l$ to denote the union of triangles in $\cT_h$ that have $Q_l$ as a common vertex. Denote the area of $\omega_l$ by $|\omega_l|$. For each vertex $Q_l$, we define the hat function $\phi_l$ so that
\begin{align*}
	\phi_l\in V_h, \ \phi_l(Q_l)=1 \mbox{ and } \phi_l(Q_m)=0 \mbox{ for } m\neq l.
\end{align*}
We have that $\phi_l$ is supported on $\omega_l$. For any function $f\in H^1(\Omega)$, its finite element approximation $f_h\in V_h$ can be written as
\begin{align*}
	f_h=\sum_{l=1}^{N_h} f(Q_l)\phi_l.
\end{align*}
We further equip $V_h$ with the inner product $(f_h,g_h)_h:V_h\times V_h\rightarrow \RR$ defined by
\begin{align*}
	(f_h,g_h)_h=\frac{1}{3}\sum_{l=1}^{N_h} |\omega_l|f_h(Q_l)g_h(Q_l), \forall f_h,g_h\in V_h.
\end{align*} 
The induced norm is defined as
\begin{align*}
	\|f_h\|_h=\sqrt{(f_h,f_h)}.
\end{align*}
Because of the eventual smoothness of solutions to the inverse iteration \eqref{eq.scheme} as shown in \cite{le2020convergence}, our mixed finite-element method uses the space $V_h$ to approximate both the solution $u$ and its second-order partial derivatives $\partial^2 u/\partial x_i\partial x_j$ for $i,j=1,2$. In the rest of this section, we denote the finite-element approximation of $v$ and $\p$ by $v_h\in V_{0h}$ and $\p_h\in (V_h)^{2\times 2}$, respectively.

\subsection{Finite element approximation of the three second-order partial derivatives}
In equation (\ref{eq.alg.2}), one needs to compute $\D^2 v^{n+1/3}$, the Hessian of $v^{n+1/3}$, which will be numerically computed, and we adopt to our current setting the \emph{double regularization} method introduced in \cite{glowinski2019finite}. 

The double regularization method is a two-step process to get a smooth approximation of $\D^2u$. In the first step, one solves
\begin{align}
	\begin{cases}
		-\varepsilon_1 \nabla^2 \pi_{ij}+\pi_{ij}=\frac{\partial^2 u}{\partial x_i\partial x_j} & \mbox{ in } \Omega,\\
		\pi_{ij}=0 &\mbox{ on } \partial \Omega,
	\end{cases}
\label{eq.double.1}
\end{align}
in which $\varepsilon_1=O(h^2)$ is a constant, $\pi_{ij}$ is a regularized approximation of $\partial^2 u/\partial x_i\partial x_j$ with zero boundary condition. Although $\pi_{ij}$ is a smooth approximation, the zero boundary condition will have a disastrous influence to the solution $u$ of our scheme, as mentioned in \cite{glowinski2019finite}. To mitigate the influence, the second step is a correction step which solves
\begin{align}
	\begin{cases}
		-\varepsilon_1 \nabla^2 D_{ij}^2u+D_{ij}^2u=\pi_{ij} & \mbox{ in } \Omega,\\
		\frac{\partial D_{ij}^2u}{\partial \n}=0 &\mbox{ on } \partial \Omega,
	\end{cases}
\label{eq.double.2}
\end{align}
where $\n$ denotes the outward normal direction of $\partial \Omega$. The resulting $D_{ij}^2u$ is the doubly regularized approximation of $\partial^2 u/\partial x_i\partial x_j$.

From the divergence theorem, one has
\begin{align}
	\begin{cases}
		\forall i,j=1,2, \ \forall v\in H^2(\Omega),\\
		\displaystyle\int_{\Omega} \frac{\partial^2 v}{\partial x_i\partial x_j} w d\x=-\frac{1}{2} \displaystyle\int_{\Omega} \left[ \frac{\partial v}{\partial x_i}\frac{\partial w}{\partial x_j}+ \frac{\partial v}{\partial x_j}\frac{\partial w}{\partial x_i}\right],\\
		\forall w\in H_0^1(\Omega).
	\end{cases}
	\label{eq.divThm}
\end{align} 
Based on equation (\ref{eq.divThm}), the discrete analogues of equations (\ref{eq.double.1})-(\ref{eq.double.2}) read as:
\begin{align}
	\begin{cases}
		\pi_{ijh}\in V_{0h},\\
		c\sum\limits_{T\in \omega_l}|T| \displaystyle\int_T \nabla \pi_{ijh} \cdot \nabla\phi_l d\x +\frac{1}{3}|\omega_l|\pi_{ijh}(Q_l)= -\frac{1}{2} \displaystyle\int_{\omega_l} \left[ \frac{\partial u_h}{\partial x_i} \frac{\partial \phi_l}{\partial x_j} +\frac{\partial u_h}{\partial x_j} \frac{\partial \phi_l}{\partial x_i} \right]d\x,\\
		\forall l=1,...,N_{0h}
	\end{cases}
\label{eq.double.1.dis}
\end{align}
and
\begin{align}
	\begin{cases}
		D_{ijh}^2 u_h\in V_h,\\
		c\sum\limits_{T\in \omega_l} |T| \displaystyle\int_{T} \nabla D_{ijh}^2 u_h\cdot \nabla\phi_l d\x +\frac{1}{3}|\omega_l|D_{ijh}^2u_h(Q_l)  = \frac{1}{3}|\omega_l|\pi_{ijh}(Q_l),\\
		\forall l=1,...,N_h,
	\end{cases}
\label{eq.double.2.dis}
\end{align}
where $c=O(1)$ is a constant.

\subsection{On the finite-element approximation of problem (\ref{eq.alg.1})}
We first rewrite equation (\ref{eq.alg.1}) in the variational form
\begin{align}
	\begin{cases}
		v^{n+1/3}\in V_{0h},\\
		\displaystyle\int_{\Omega} v^{n+1/3}\psi d\x+\tau \displaystyle\int_{\Omega} (\varepsilon\I+\cof(\p^n))\nabla v^{n+1/3} \cdot \nabla \psi d\x =2\displaystyle\int_{\Omega} f^k\psi d\x,\\
		\forall \psi\in V_{0h}.
	\end{cases}
	\label{eq.ellip.var}
\end{align}
If $\p^n$ is semi--positive definite, then problem (\ref{eq.ellip.var}) admits a unique solution. Denote $\M=\varepsilon\I+\cof(\p^n_h)$. The discrete analogue of equation (\ref{eq.ellip.var}) reads as
\begin{align}
	\begin{cases}
		v_h^{n+1/3}\in V_{0h},\\
		\frac{1}{3}|\omega_l|v_h^{n+1/3}(Q_l) +\tau \sum\limits_{m=1}\limits^{N_{0h}} \left(v_h^{n+1/3}(Q_m)\displaystyle\int_{\omega_l\cap \omega_m} \M\nabla \phi_m\cdot \nabla \phi_l dx \right)= \frac{2}{3}|\omega_l|f^k(Q_l),\\
		\forall l=1,...,N_{0h}.
	\end{cases}
\label{eq.ellip}
\end{align}
Solving problem (\ref{eq.ellip}) is equivalent to solving a sparse linear system, for which many efficient solvers, such as the Cholesky decomposition, can be used.

\subsection{On the finite element approximation of problem (\ref{eq.alg.2})}
\label{sec.projP}
We first define the projection operator $P_+$ that projects $2\times 2$ real symmetric matrices to the set of real symmetric semi-positive definite matrices. Let $\A$ be a $2\times 2$ real symmetric matrix. By spectral decomposition, there exists a $2\times 2$ orthogonal matrix $\bS$ so that $\A=\bS\bm{\Lambda}\bS^{-1}$, where
\begin{align*}
	\bm{\Lambda}=\begin{bmatrix}
		\lambda_1 &0\\
		0 & \lambda_2
	\end{bmatrix}
\end{align*}
with $\lambda_1,\lambda_2$ being eigenvalues of $\A$. If $\A$ is semi--positive definite, one has $\lambda_1,\lambda_2\geq 0$. Therefore we define $P_+$ as
\begin{align*}
	P_+(\A)=\bS\begin{bmatrix}
		\max(\lambda_1,0) &0 \\
		0 & \max(\lambda_2,0)
	\end{bmatrix}\bS^{-1}.
\end{align*}

In equation (\ref{eq.alg.2}), we compute 
\begin{align}
	\p_h^{n+1}=P_+\left(e^{-\gamma\tau}\p_h^{n}+(1-e^{-\gamma\tau})
	\begin{bmatrix}
		D^2_{11h}v_h^{n+1/3} & D^2_{12h}v_h^{n+1/3}\\
		D^2_{21h}v_h^{n+1/3} & D^2_{22h}v_h^{n+1/3}
	\end{bmatrix}\right),
\end{align}
where the entries $D^2_{ijh}v_h^{n+1/3}$ are computed using equations (\ref{eq.double.1.dis})-(\ref{eq.double.2.dis}).

\subsection{On the finite element approximation of problem (\ref{eq.alg.3})}
According to equation (\ref{eq.v.proj}), we compute $v_h^{n+1}$ as
\begin{align}
	v_h^{n+1}=\frac{v_h^{n+1/3}}{\left(\sum\limits_{l=1}\limits^{N_{0h}} \frac{1}{3}|\omega_l|\left(v_h^{n+1/3}(Q_l)\right)^2\right)^{1/2}}.
\end{align}

\subsection{On the finite element approximation of equation (\ref{eq.Rayleigh.div})}
For any $u_h\in H^1_0(\Omega)$, the discrete analogue of equation (\ref{eq.Rayleigh.div}) reads as
\begin{align}
	R(u_h)=-\frac{ \sum\limits_{m,l=1}\limits^{N_{0h}} u_h(Q_m)u_h(Q_l)\displaystyle\int_{\omega_l\cap \omega_m} (\cof(\D^2_h u_h(Q_m)) \nabla \phi_m)\cdot \nabla \phi_l d\x}{\sum\limits_{m=1}\limits^{N_{0h}}\frac{2}{3}|\omega_m|(-u_h(Q_m))^3},
	\label{eq.Rayleigh.div.dis}
\end{align}
where $\D^2_h u_h$ is the finite-element approximation of $\D^2u$ computed using equations (\ref{eq.double.1.dis}) and (\ref{eq.double.2.dis}).

Note that if $u$ is an eigenfunction of the Monge--Amp\`{e}re equation (\ref{eq.MAEV}), by Theorem \ref{thm.eigen}, one can compute the eigenvalue as $\lambda_{\rm MA}=\inf_{u\in \cK} R(u)$. Therefore, for every time step, we can compute the approximate `eigenvalue' corresponding to $u_h^k$ as
\begin{align*}
	\lambda_h^k=R(u_h^k)
\end{align*}
and monitor the evolution of $\lambda_h^k$, which will monotonically converge to $\lambda_{\rm MA}$ as shown in \cite{le2020convergence}. 

\subsection{On the finite element approximation of the initial condition}
Denote the finite element of $u_0$ and $(v_0,\p_0)$ by $u_{0h}$ and $(v_{0h},\p_{0h})$, respectively. The discrete analogue of the initial condition  (\ref{eq.initial.outer}) reads as
\begin{align*}
	\begin{cases}
		u_{0h}\in V_{0h},\\
		\sum\limits_{m=1}\limits^{N_{0h}} u_{0h}(Q_m)\displaystyle\int_{\omega_l\cap \omega_m} \nabla \phi_m \cdot \nabla \phi_l d\x =\frac{\eta}{3}|\omega_l|,\\
		\forall l=1,...,N_{0h}.
	\end{cases}
\end{align*}

For $(v_{0h},\p_{0h})$, we set
\begin{align*}
	v_{0h}=u_h^k,\quad \p_{0h}=\D^2_h v_{0h},
\end{align*}
where $\D^2_h$ is the double regularization approximation using equations (\ref{eq.double.1.dis})-(\ref{eq.double.2.dis}).

\section{Numerical experiments}
\label{sec.experiments}

We demonstrate the efficiency of scheme (\ref{eq.alg.1})-(\ref{eq.alg.3}) by several numerical experiments. We set the stopping criterion as 
\begin{align}
	\|u_h^{k+1}-u_h^{k}\|_h< \xi
\end{align}
for some small $\xi>0$. Without specification, in all of our experiments, we set $\xi=10^{-6}$, $\varepsilon=2h^2$, and $c=2$, where $\varepsilon$ and $c$ are regularization parameters in equation (\ref{eq.cMA.2}) and scheme (\ref{eq.double.1.dis})-(\ref{eq.double.2.dis}), respectively.

When the exact solution, denoted by $u_h^*$, is given, we define the $L^2$ error and $L^{\infty}$ error of $u_h$ as
\begin{align}
	\|u_h-u_h^*\|_h \quad \mbox{ and } \quad \max_m |u_h(Q_m)-u_h^*(Q_m)|,
\end{align}
respectively.

Algorithm \ref{alg.1} consists of two iterations: the outer iteration for $u$ and the inner iteration for $v$ and $\p$. Since both $u$ and $v$ are estimates of the solution of equation (\ref{eq.MAEV}), it is not necessary to solve every inner iteration until steady state. Instead, one can just solve the inner iteration for a few steps. In our experiments, we observe that just 1 iteration step for the inner iteration is sufficient for our algorithm to converge. Thus in all of our experiments, we solve the inner iteration for only 1 step in each outer iteration.

\begin{figure}
	\begin{tabular}{cc}
		(a) & (b)\\
		\includegraphics[width=0.25\textwidth]{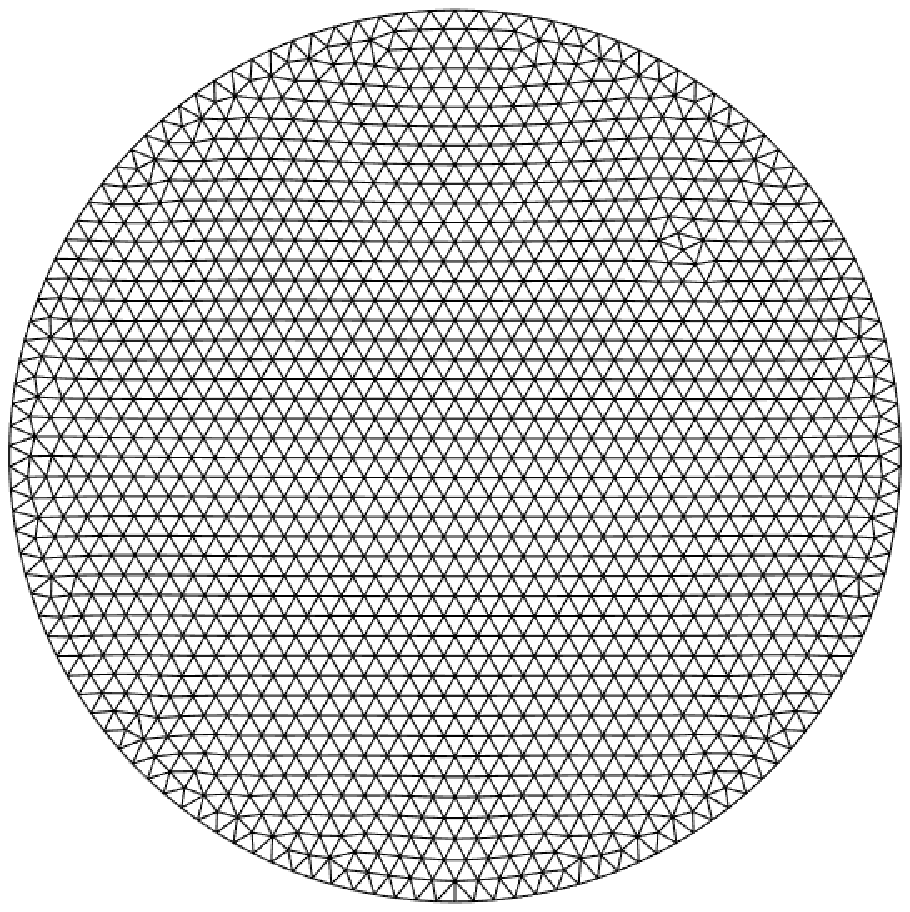}&
		\includegraphics[width=0.25\textwidth]{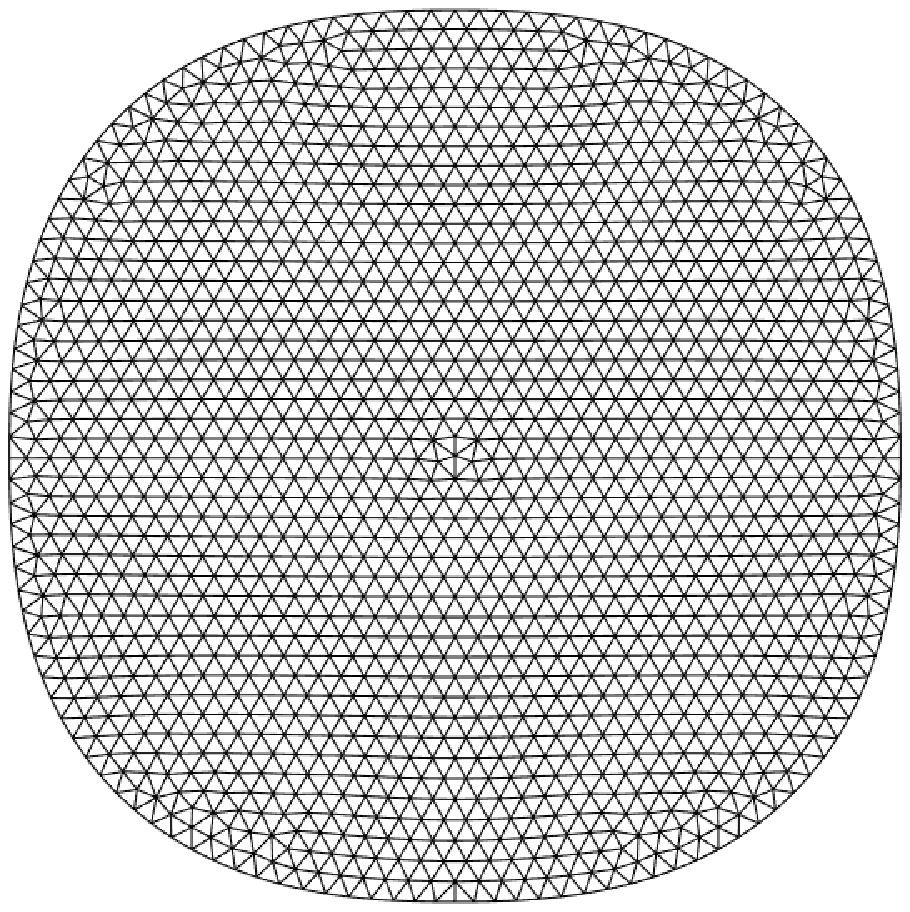}\\
		(c) & (d)\\
		\includegraphics[width=0.3\textwidth]{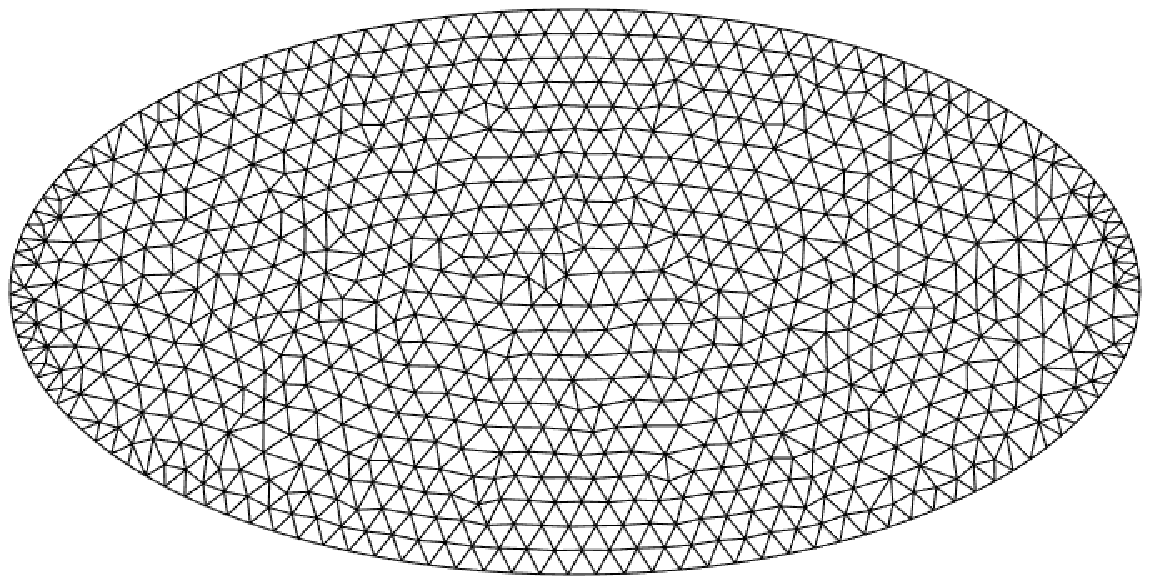}&
		\includegraphics[width=0.3\textwidth]{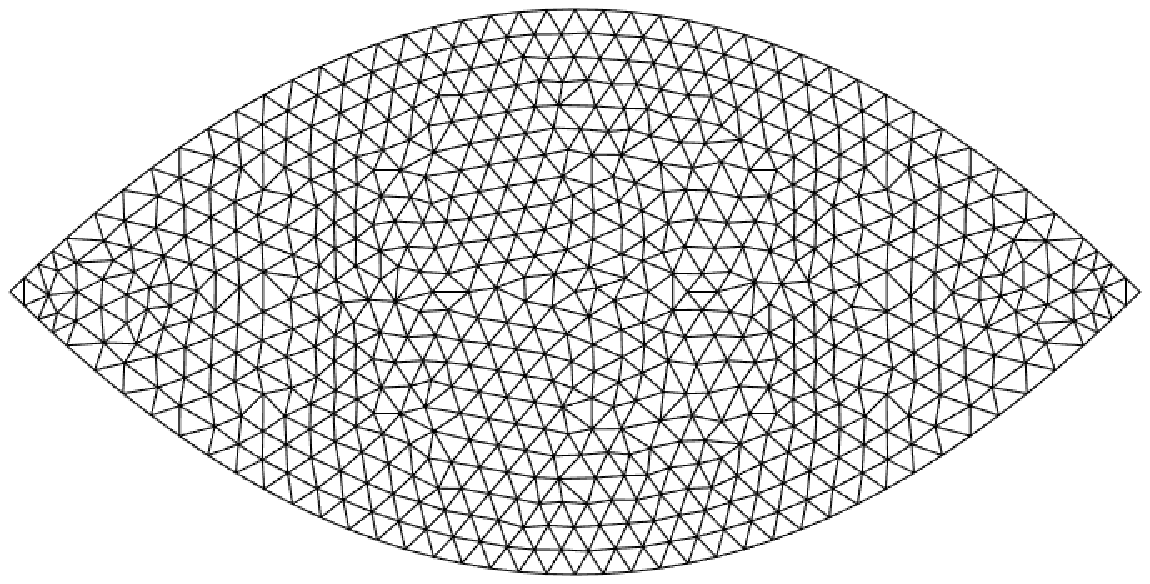}
	\end{tabular}
	\caption{The triangulation of domains used in the examples. (a) The unit disk domain (\ref{eq.disk}) with $h=1/20$. (b) The smoothed square domain (\ref{eq.sq}) with $h=1/20$. (c) The ellipse domain (\ref{eq.ellipse}) with $h=1/20$. (d) The eye-shape domain (\ref{eq.eye}) with $h=1/40$.}
	\label{fig.domain}
\end{figure}

\subsection{Example 1}
In the first example, we test our algorithm on the unit disk
\begin{align}
	\Omega=\{(x_1,x_2): x_1^2+x_2^2 < 1\}.
	\label{eq.disk}
\end{align}
The triangulation of the domain with $h=1/20$ is visualized in Figure \ref{fig.domain}(a).

\begin{figure}
	\begin{tabular}{ccc}
		(a) & (b) & (c)\\
		\includegraphics[width=0.3\textwidth]{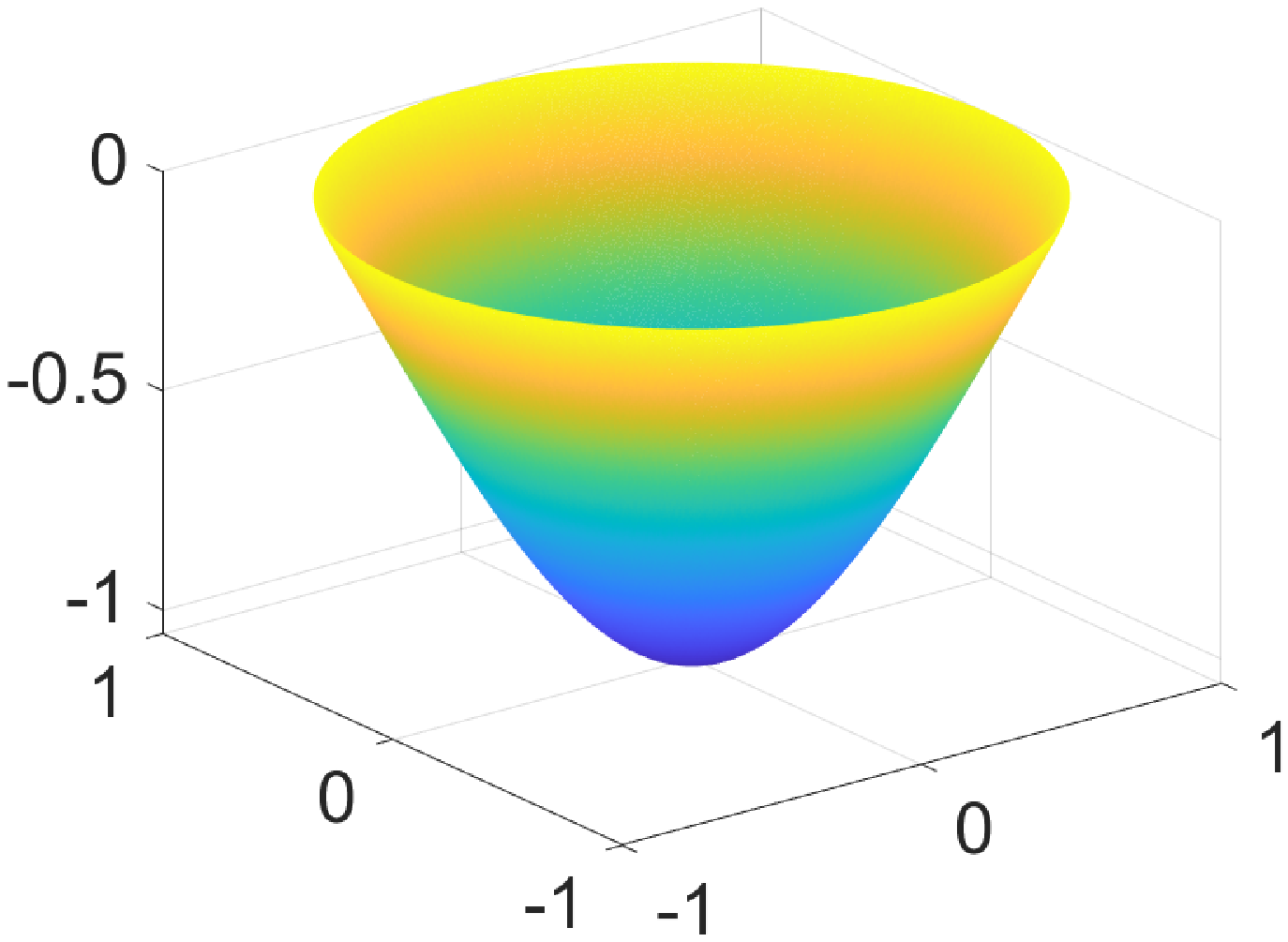}&
		\includegraphics[width=0.3\textwidth]{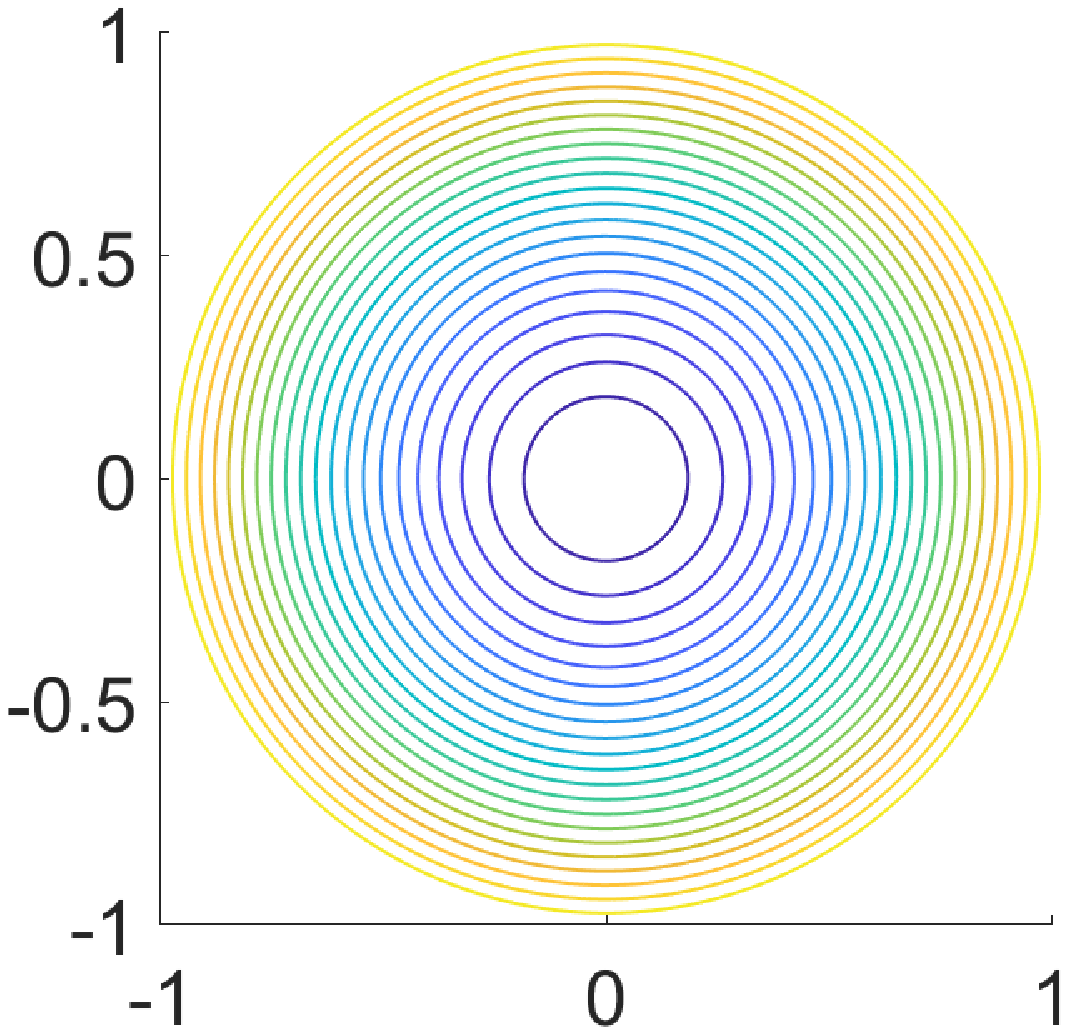}&
		\includegraphics[width=0.3\textwidth]{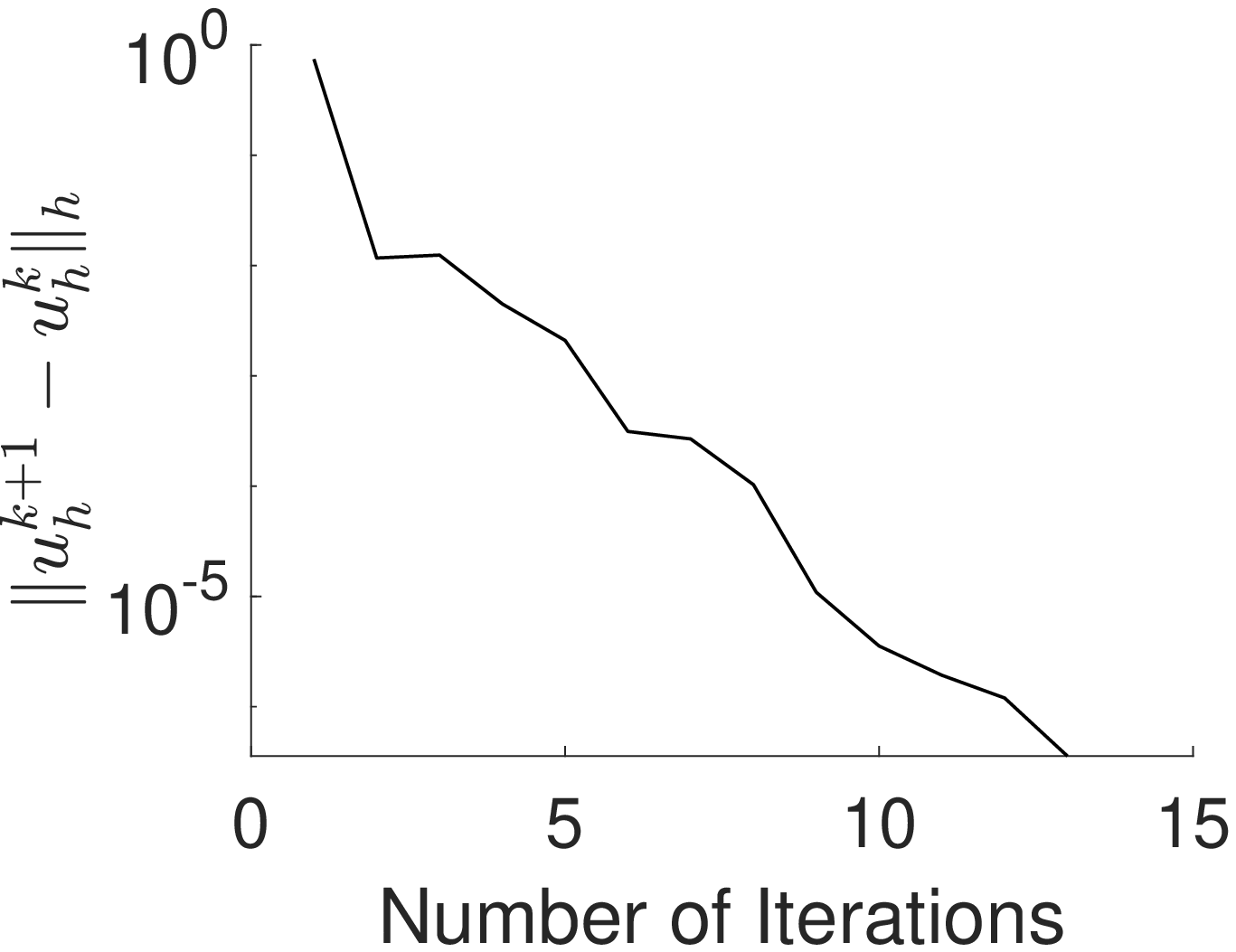} \\
		(d) & (e) & (f)\\
		\includegraphics[width=0.3\textwidth]{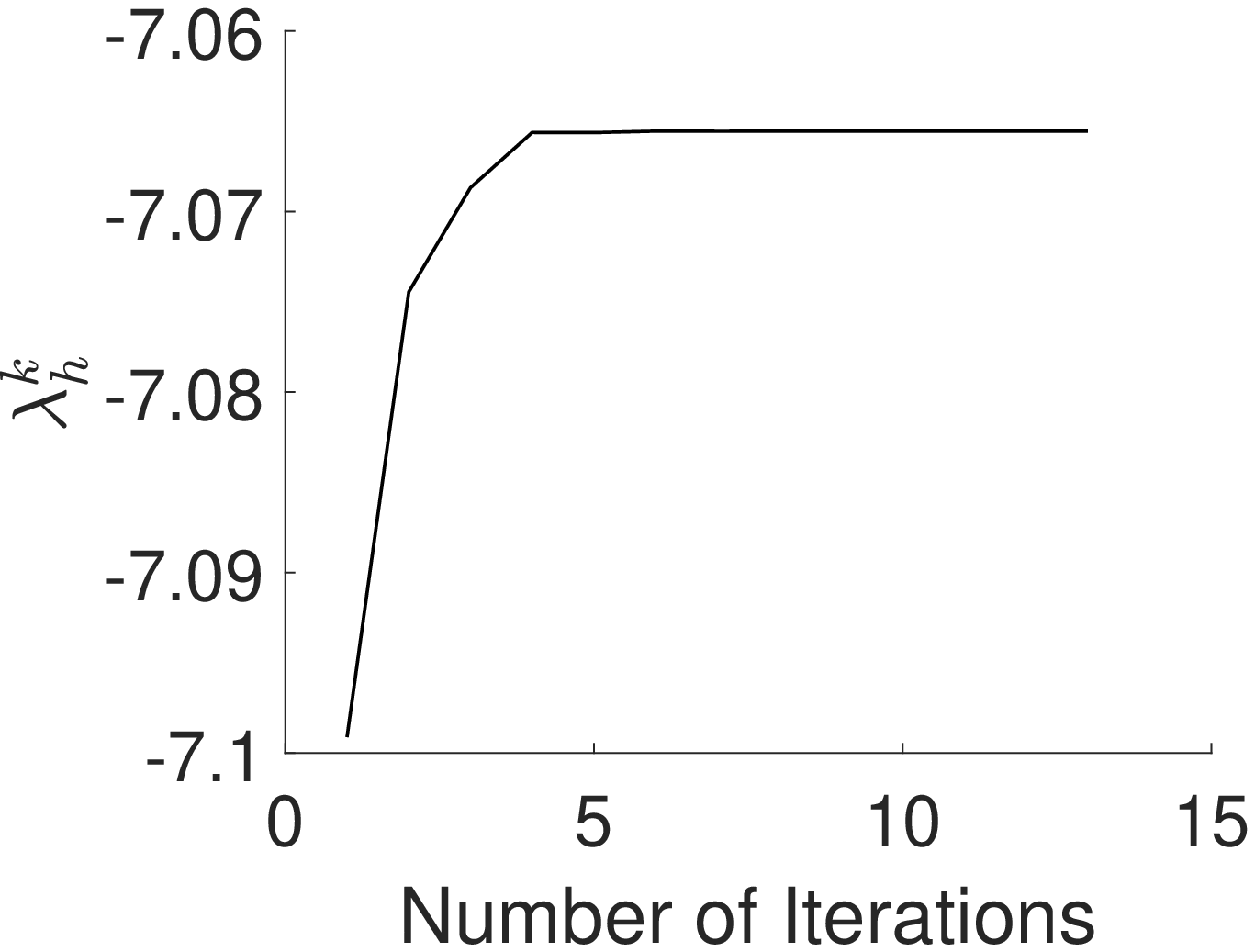}&
		\includegraphics[width=0.3\textwidth]{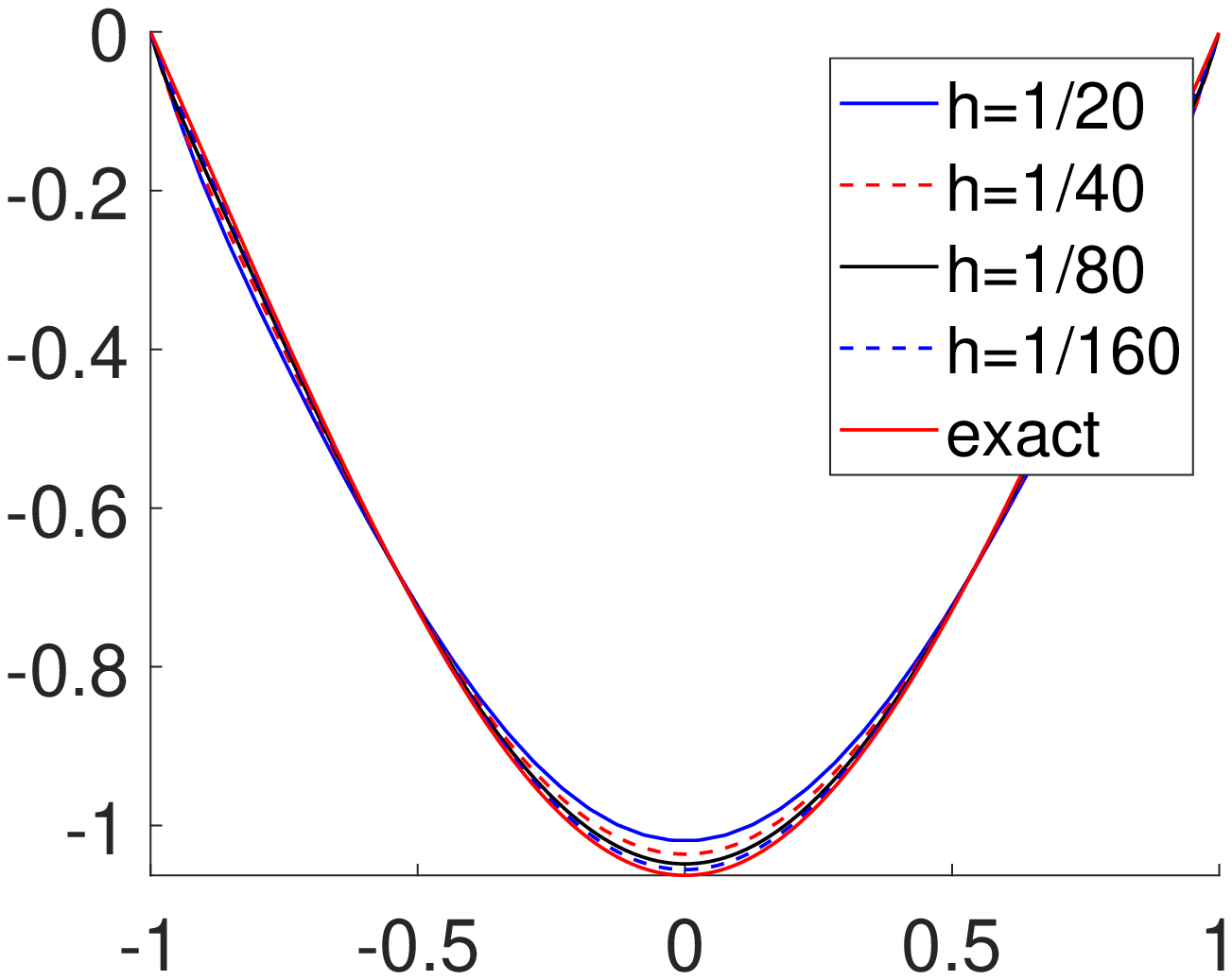}&
		\includegraphics[width=0.3\textwidth]{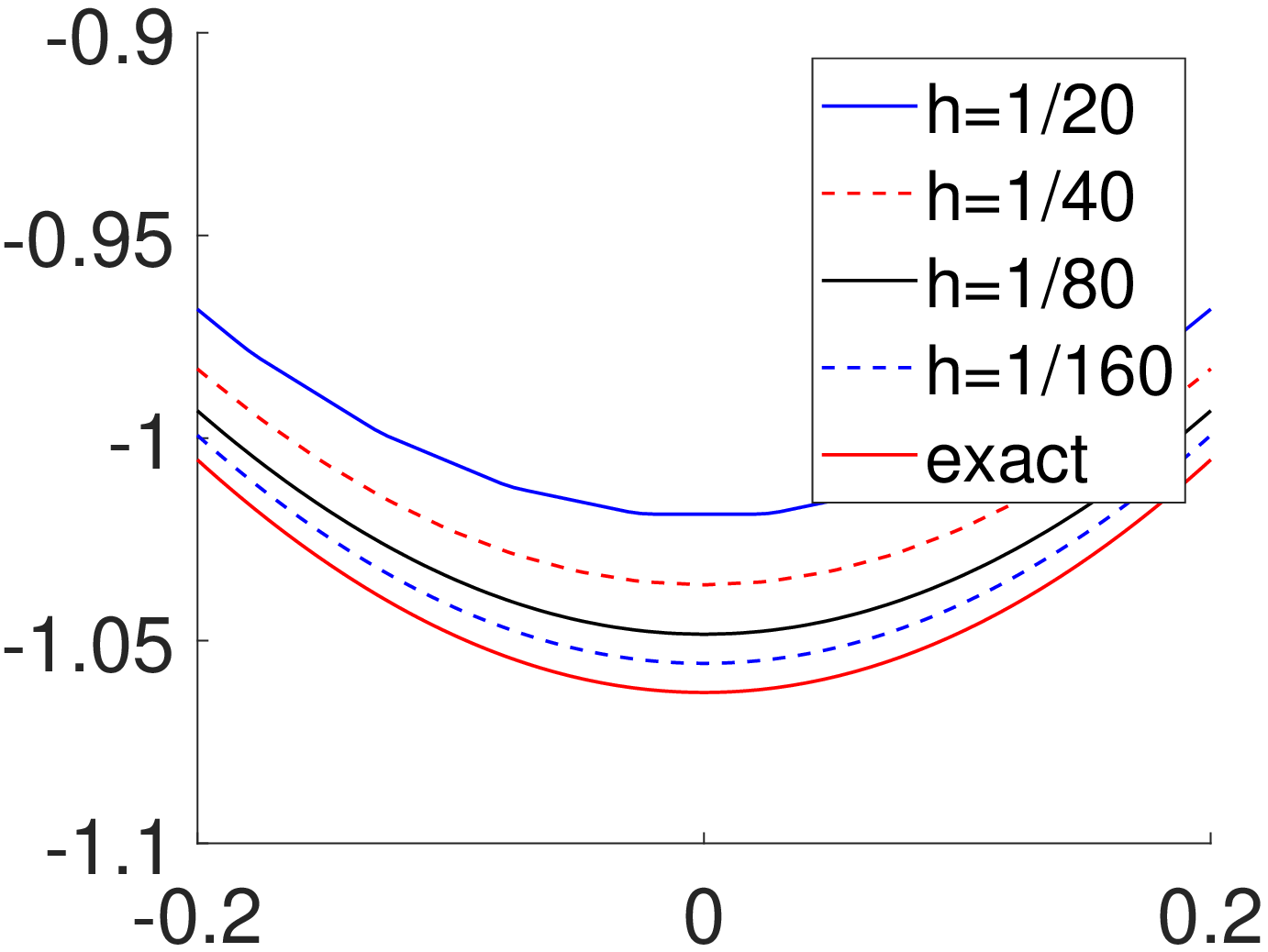}
	\end{tabular}
	\caption{The unit disk domain (\ref{eq.disk}). (a) The computed result with $h=1/80$. (b) The contour of (a). (c) The history of the error $\|u_h^{k+1}-u_h^k\|_h$ with $h=1/80$. (d) The history of the computed eigenvalue $\lambda_h^k$ with $h=1/80$. (e) Comparison of the cross sections along $x_2=0$ of the computed solution with various $h$. (f) Zoomed plot of the bottom region of (e).}
	\label{fig.disk}
\end{figure}

In this case, equation (\ref{eq.MAEV}) has a radial solution. Let $r=\sqrt{x_1^2+x_2^2}$. For a radial function $g(r)$, one has $\det \D^2 g=\frac{g'g''}{r}$. Therefore, we write the solution to equation (\ref{eq.MAEV}) as $u(r)$, which satisfies
\begin{align}
	\begin{cases}
		u\leq 0, \ \lambda> 0,\\
		u'\,u''=-\lambda\, r\,u^2\; \mbox{ in } (0,1),\\
		u'(0)=0,\ u(1)=0,\\
		2\pi \displaystyle\int_0^1 |u|^2\,r\,dr=1.
	\end{cases}
\label{eq.ode}
\end{align}
Using a shooting method, we can solve the ODE problem (\ref{eq.ode}) very accurately. The `exact' solution verifies $u(\zero)\approx-1.0628$ and $\lambda\approx 7.4897$. On the domain (\ref{eq.disk}), we test our algorithm with $h=1/20,1/40,1/80$ and $1/160$. In Figure \ref{fig.disk}(a)--(d), we show results with $h=1/80$. Our numerical result is visualized in Figure \ref{fig.disk}(a). The contour of Figure \ref{fig.disk}(a) is shown in Figure \ref{fig.disk}(b). Our result is a smooth radial function, whose contour consists of several circles with the same center. The convergence histories of the error $\|u_h^{k+1}-u_h^k\|_h$ and the computed eigenvalue are shown in Figure \ref{fig.disk}(c) and Figure \ref{fig.disk}(d), respectively. Linear convergence is observed for the error $\|u_h^{k+1}-u_h^k\|_h$, and the convergence rate is approximately 0.47. The computed eigenvalue converges with just 5 iterations. In Figure \ref{fig.disk}(e), we show the cross sections of the results with various $h$ along $x_2=0$. As $h$ goes to 0, our computed solution converges to the exact solution. For better visualization, the zoomed bottom region of Figure \ref{fig.disk}(e) is shown in Figure \ref{fig.disk}(f).

To quantify the convergence of the proposed algorithm, we present in Table \ref{tab.disk} the number of iterations needed for convergence, $L^2$- and $L^{\infty}$-errors, computed eigenvalues and the minimal value of the computed solution with various $h$. For all resolutions of mesh, 13 iterations are sufficient for the algorithm to converge. As $h$ goes to zero, the convergence rate of the $L^2$- and $L^{\infty}$-error goes to 1, and the computed eigenvalue and the minimal value converge to the exact solutions. The eigenvalue $\lambda_h$ converges linearly to the exact eigenvalue with an error of $O(h)$.

We next compare Algorithm \ref{alg.1} with the method proposed in \cite{glowinski2020numerical}. For the method from \cite{glowinski2020numerical}, we have to use small time steps to make sure that the method does converge. In the numerical experiment, we set the time step as $h/2$ and stopping criterion as $10^{-6}$. Note that the method from \cite{glowinski2020numerical} finds the solution of equation (\ref{eq.MAEV}) with $\|u_h\|_3=1$. When computing the $L^2$- and $L^{\infty}$- errors, we first normalize the solution so that $\|u_h\|_2=1$ and we then compute the errors. The comparisons are shown in Table \ref{tab.disk.compare}. For both $L^2$- and $L^{\infty}$- errors, both algorithms have errors with similar magnitudes. 
We compare the computational efficiency between the two algorithms in Table \ref{tab.disk.compare.cpu}. The number of iterations used by Algorithm \ref{alg.1} is independent of the mesh resolution, while the number of iterations used by \cite{glowinski2020numerical} grows approximately linearly with $1/h$. For the CPU time, Algorithm (\ref{alg.1}) is also much faster than the method in \cite{glowinski2020numerical}. Note that in Algorithm (\ref{alg.1}), the constraint $\|u_h\|_2=1$ is enforced by the projection step (\ref{eq.alg.3}). In \cite{glowinski2020numerical}, the constraint is $\|u_h\|_3=1$, which was enforced by a sequential quadratic programming algorithm, which in turn uses around 15 iterations in each outer iteration.

\begin{table}[t!]
	\centering
	\begin{tabular}{|c|c|c|c|c|c|c|c|c|}
		\hline
		$h$& \# Iter. & $\|u_h^{k+1}-u_h^k\|_{h}$ & $L^2$-error & rate &$L^{\infty}$-error &rate&$\lambda_h$& $\min u_h$ \\
		\hline
		1/20 & 13 & 2.13$\times 10^{-7}$ &$4.91\times10^{-2}$ & & $4.29\times10^{-2}$ & & 5.9716 & -1.0189\\
		\hline
		1/40 & 13 & 2.91$\times 10^{-7}$ & $3.36\times10^{-2}$ & 0.54& $3.04\times10^{-2}$& 0.50 &6.6656 & -1.0362\\
		\hline
		1/80 & 13 & 3.56$\times 10^{-7}$ & $1.94\times10^{-2}$ & 0.79& $1.86\times10^{-2}$ & 0.71&7.0655 & -1.0484\\
		\hline
		1/160 & 13 & 4.04$\times 10^{-7}$ & $1.01\times10^{-2}$ & 0.94& $1.03\times10^{-2}$ & 0.85&7.2816 & -1.0556\\
		\hline
	\end{tabular}
	\caption{The unit disk domain (\ref{eq.disk}). Variations with $h$ of the number of iterations necessary to achieve convergence (2nd column), of the $L^2$ and $L^{\infty}$ approximation errors and of the associated convergence rates (columns 4, 5, 6 and 7), of the computed eigenvalue (8th column) and of the minimal value of $u_h$ over $\Omega$ (that is $u_h(\mathbf{0})$) (9th column). The exact eigenvalue is around $7.4897$. The minimal value of the exact solution is  around $-1.0628$.}
	\label{tab.disk}
\end{table}

\begin{table}[t!]
	\centering
	\begin{tabular}{|c|c|c|c|c||c|c|c|c|}
		\hline
		& \multicolumn{4}{c||}{Algorithm (\ref{alg.1})} & \multicolumn{4}{c|}{Method from \cite{glowinski2020numerical}}\\
	\hline
		$h$& $L^2$-error  & rate &$L^{\infty}$-error &rate& $L^2$-error  & rate &$L^{\infty}$-error &rate \\
		\hline
		1/20 &  $4.91\times10^{-2}$ & & $4.29\times10^{-2}$ & & $4.01\times10^{-2}$ & & $8.40\times10^{-2}$ & \\
		\hline
		1/40 &   $3.36\times10^{-2}$ & 0.54& $3.04\times10^{-2}$& 0.50 &$2.33\times10^{-2}$ & 0.78& $4.00\times10^{-2}$ & 1.07\\
		\hline
		1/80 &   $1.94\times10^{-2}$ & 0.79& $1.86\times10^{-2}$ & 0.71&$1.37\times10^{-2}$ & 0.76& $2.05\times10^{-2}$ & 0.96\\
		\hline
		1/160 &   $1.01\times10^{-2}$ & 0.94& $1.03\times10^{-2}$ & 0.85&$7.55\times10^{-3}$ & 0.86& $1.08\times10^{-2}$ & 0.92\\
		\hline
	\end{tabular}
	\caption{The unit disk domain (\ref{eq.disk}). Variations with $h$ of the number of iterations necessary to achieve convergence (2nd column), of the $L^2$ and $L^{\infty}$ approximation errors and of the associated convergence rates (columns 4, 5, 6 and 7), of the computed eigenvalue (8th column) and of the minimal value of $u_h$ over $\Omega$ (that is $u_h(\mathbf{0})$) (9th column). The exact eigenvalue is around $7.4897$. The minimal value of the exact solution is  around $-1.0628$.}
	\label{tab.disk.compare}
\end{table}

\begin{table}[t!]
	\centering
	\begin{tabular}{|c|c|c||c|c|}
		\hline
		& \multicolumn{2}{c||}{Algorithm (\ref{alg.1})} & \multicolumn{2}{c|}{Method from \cite{glowinski2020numerical}}\\
		\hline
		$h$& \# Iter.   & CPU time &\# Iter.   & CPU time \\
		\hline
		1/20 &  13 & 1.44 & 62 & 3.55\\
		\hline
		1/40 &  13 & 4.58 & 101 & 22.39\\
		\hline
		1/80 &  13 & 18.35 & 151 & 138.47 \\
		\hline
		1/160 &  13 & 83.95 & 263 & 1206.96\\
		\hline
	\end{tabular}
	\caption{The unit disk domain (\ref{eq.disk}). Comparison of the number of iterations and the CPU time needed by Algorithm \ref{alg.1} and the method in \cite{glowinski2020numerical} for convergence.}
	\label{tab.disk.compare.cpu}
\end{table}

\subsection{Example 2}
\begin{figure}
	\begin{tabular}{ccc}
		(a) & (b) & (c)\\
		\includegraphics[width=0.3\textwidth]{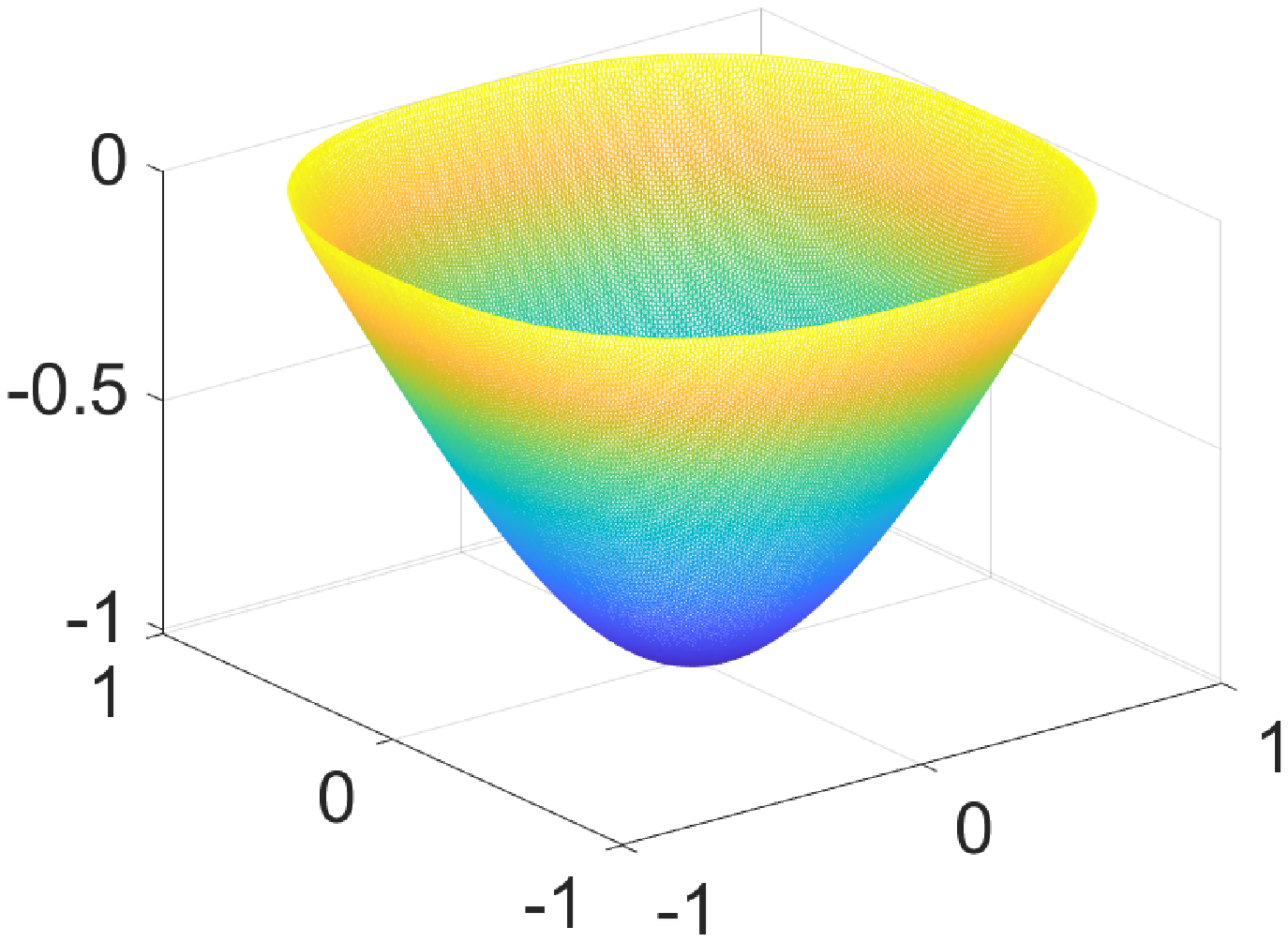}&
		\includegraphics[width=0.3\textwidth]{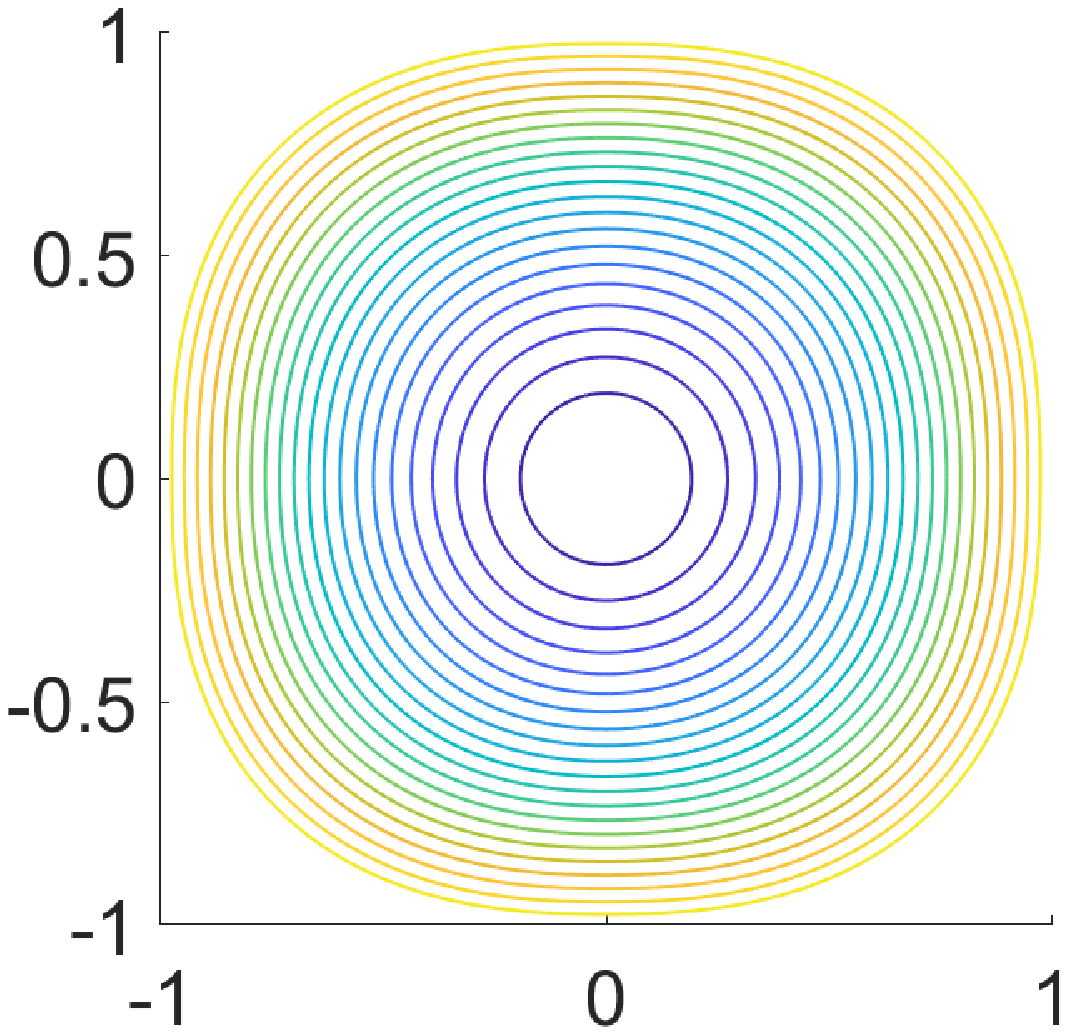}&
		\includegraphics[width=0.3\textwidth]{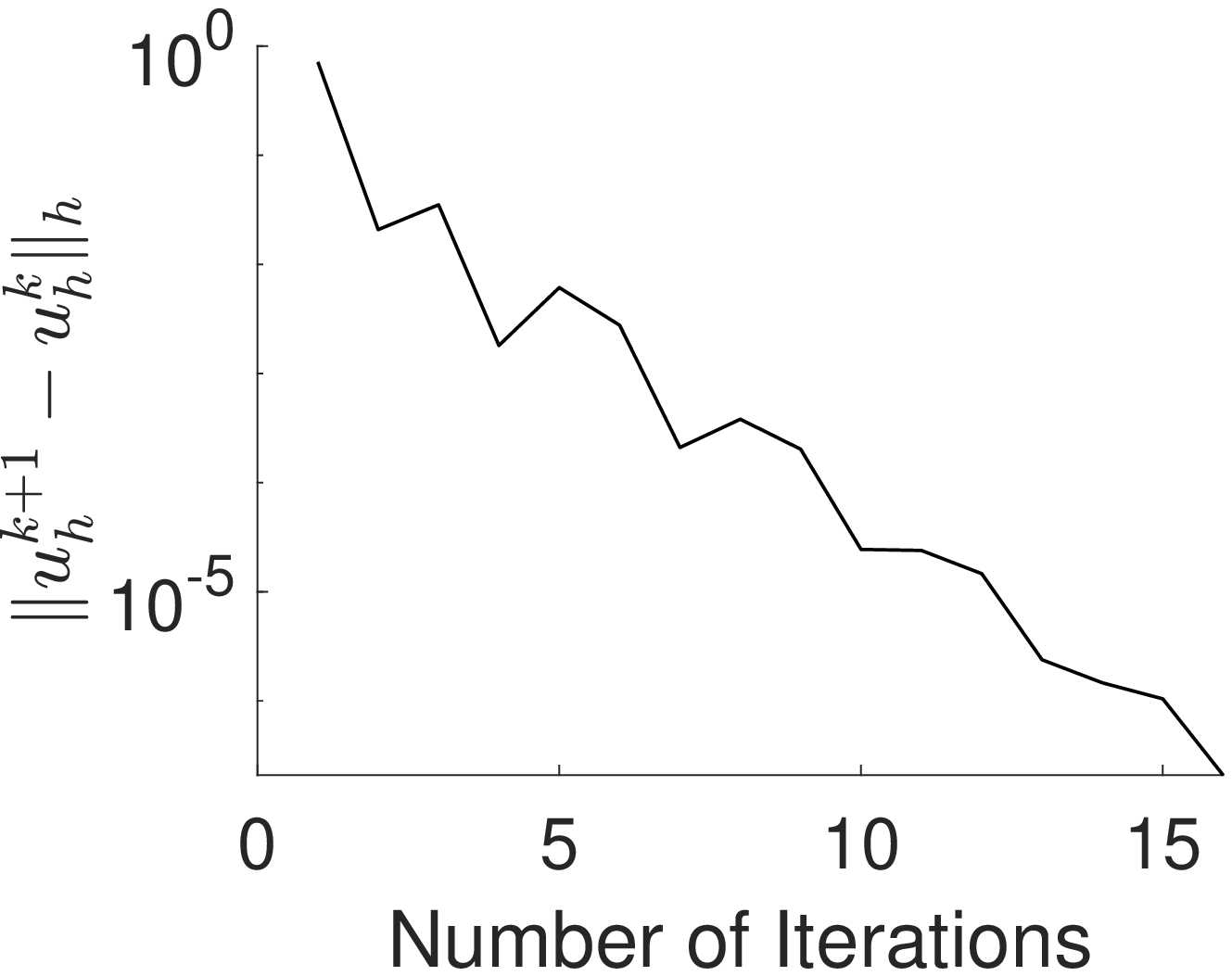} \\
		(d) & (e) & (f)\\
		\includegraphics[width=0.3\textwidth]{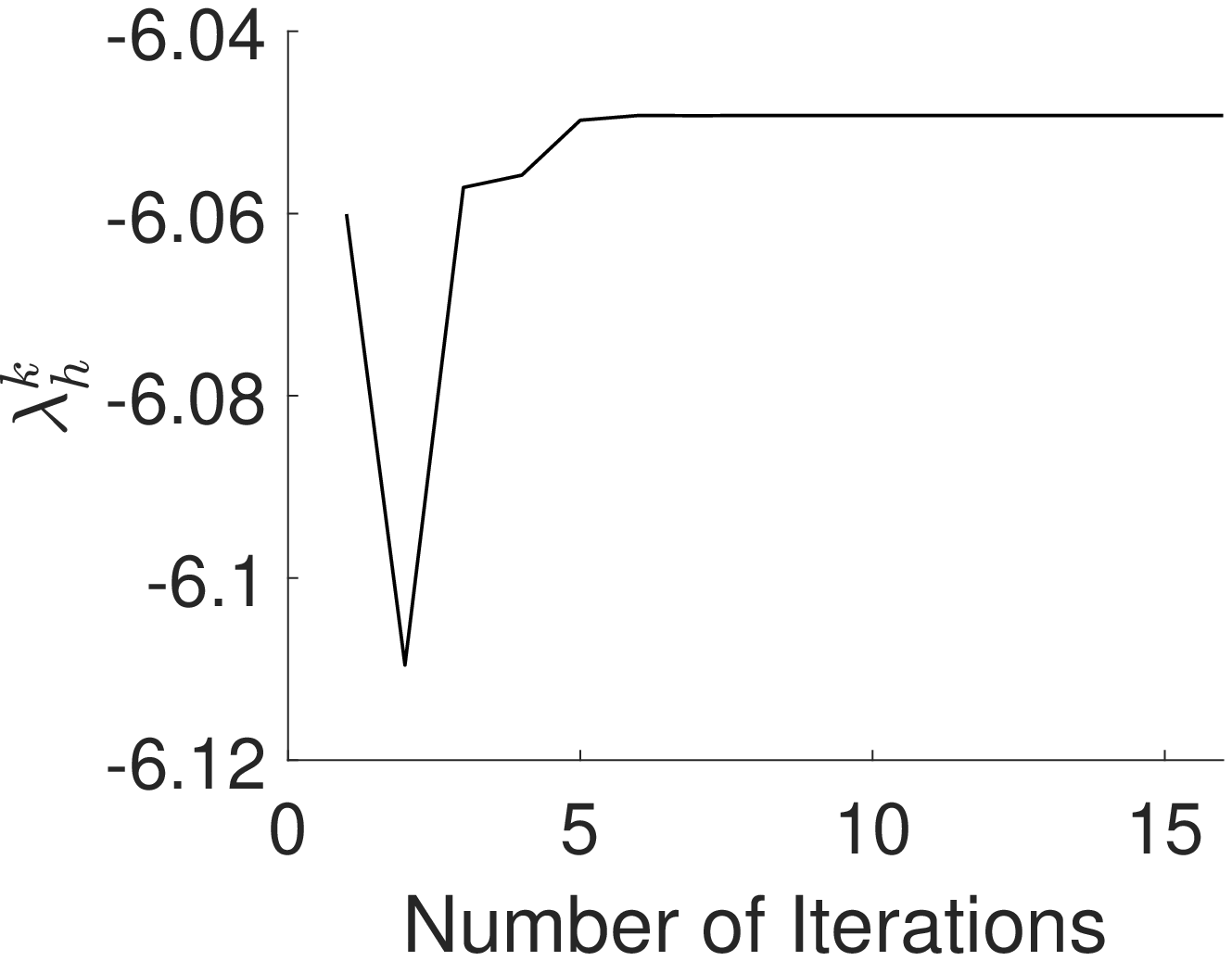}&
		\includegraphics[width=0.3\textwidth]{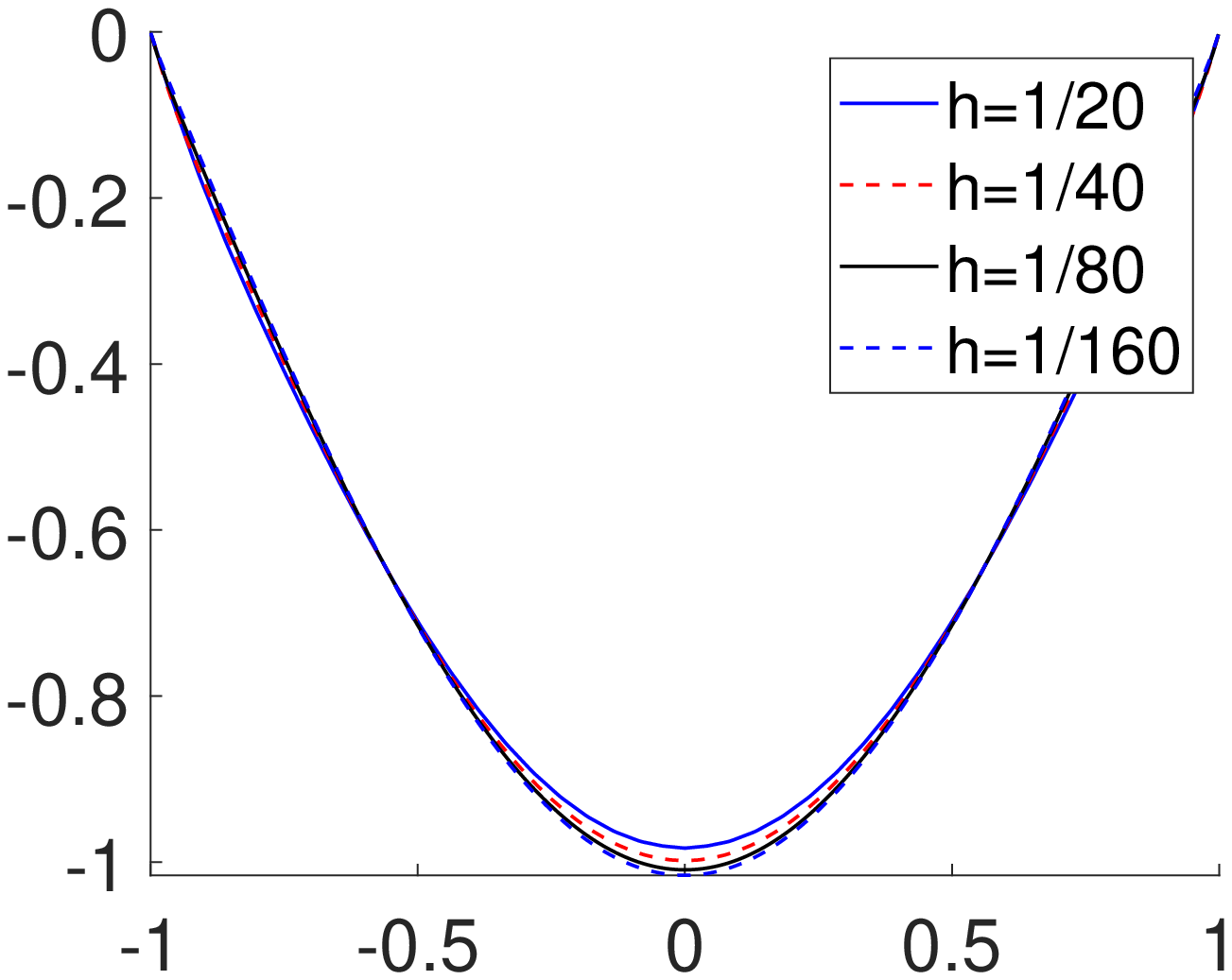}&
		\includegraphics[width=0.3\textwidth]{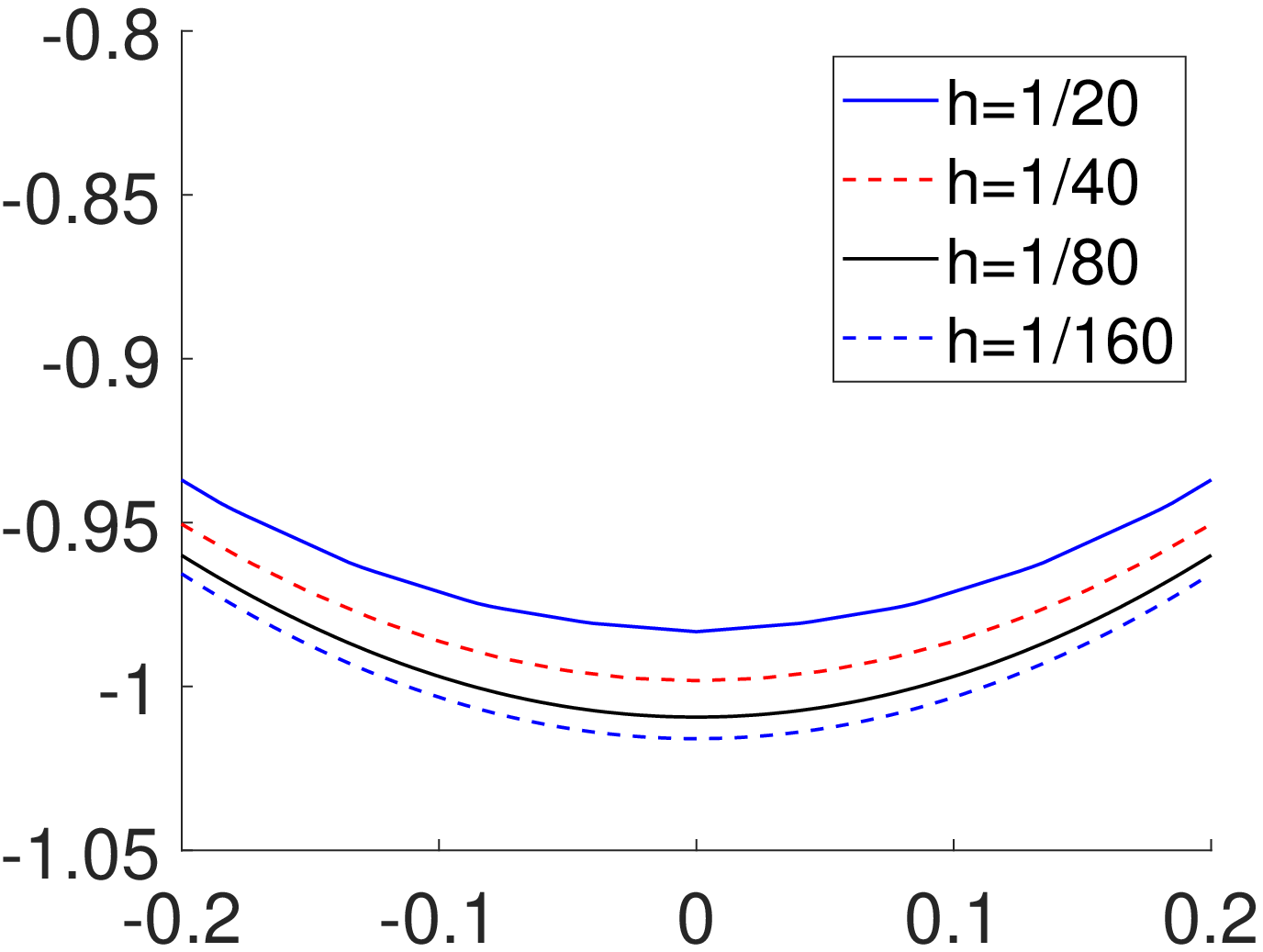}
	\end{tabular}
	\caption{The smoothed square domain (\ref{eq.sq}). (a) The computed result with $h=1/80$. (b) The contour of (a). (c) The convergence history of the errors $\|u_h^{k+1}-u_h^k\|_h$. (d) The history of the computed eigenvalue $\lambda_h^k$. (e) Comparison of the cross sections along $x_2=0$ of the computed solution with various $h$. (f) Zoomed plot of the bottom region of (e).}
	\label{fig.sq}
\end{figure}

In the second example, we consider the convex smoothed square domain
\begin{align}
	\Omega=\left\{(x_1,x_2): |x_1|^{2.5}+|x_2|^{2.5}< 1\right\}.
	\label{eq.sq}
\end{align}
The triangulation of the domain with $h=1/20$ is visualized in Figure \ref{fig.domain}(b), which has a shape between the unit disk and a square. We test our algorithm with $h$ varying from $h=1/20$ to $h=1/160$. Similar to our settings in the previous example, we set the stopping criterion $\xi=10^{-6}$. The time step is set as $\tau=1/2$. Our results with $h=1/80$ are visualized in Figure \ref{fig.sq}(a)--(d). Our computed  solution is shown in Figure \ref{fig.sq}(a), whose contour is shown in Figure \ref{fig.sq}(b). Again, our solution is very smooth. The convergence histories of the error $\|u_h^{k+1}-u_h^k\|_h$ and the computed eigenvalues $\lambda_h^k$ are shown in Figure \ref{fig.sq}(c) and Figure \ref{fig.sq}(d), respectively. The error $\|u_h^{k+1}-u_h^k\|_h$ converges linearly with a rate of 0.38. In this numerical experiment, the stopping criterion is satisfied after 16 iterations. The computed eigenvalue achieves its steady state with about 5 iterations. With various $h$, the comparison of cross sections of our results along $x_2=0$ is shown in Figure \ref{fig.sq}(e)--(f). As $h$ goes to 0, the convergence of the solution along cross sections is observed. 

We then report the computational cost and convergence behavior of the computed eigenvalue and minimal value with various $h$ in Table \ref{tab.sq}. The convergence of the eigenvalue is similar to that in \cite{glowinski2020numerical}: the eigenvalue $\lambda_h$ converges to $\lambda$ uniformly in the rate $\lambda_h\approx \lambda-ch$ with $\lambda\approx 6.4, c\approx26$. In terms of the computational cost, Algorithm \ref{alg.1} is very efficient since all experiments used less than 20 iterations to satisfy the stopping criterion.

\begin{table}[t!]
	\centering
	\begin{tabular}{|c|c|c|c|c|}
		\hline
		$h$& \# Iter.  & $\|u_h^{k+1}-u_h^k\|_{h}$ & $\lambda_h$& $\min u_h$\\
		\hline
		1/20 & 14 &6.05$\times 10^{-7}$ & 5.17 & -0.9833\\
		\hline
		1/40 & 14 & 8.00$\times 10^{-7}$ & 5.72 & -0.9982\\
		\hline
		1/80 & 16 & 2.08$\times 10^{-7}$ & 6.05 & -1.0094\\
		\hline
		1/160 & 18 & 7.77$\times 10^{-7}$ & 6.22 & -1.0159\\
		\hline
	\end{tabular}
	\caption{The smoothed square domain (\ref{eq.sq}). Variations with $h$ of the number of iterations necessary to achieve convergence (2nd column), of the computed eigenvalue (4th column) and of the minimal value of $u_h$ over $\Omega$ (that is $u_h(\mathbf{0})$) (5th column).}
	\label{tab.sq}
\end{table}

\subsection{Example 3}
\begin{figure}
	\begin{tabular}{ccc}
		(a) & (b) & (c)\\
		\includegraphics[width=0.3\textwidth]{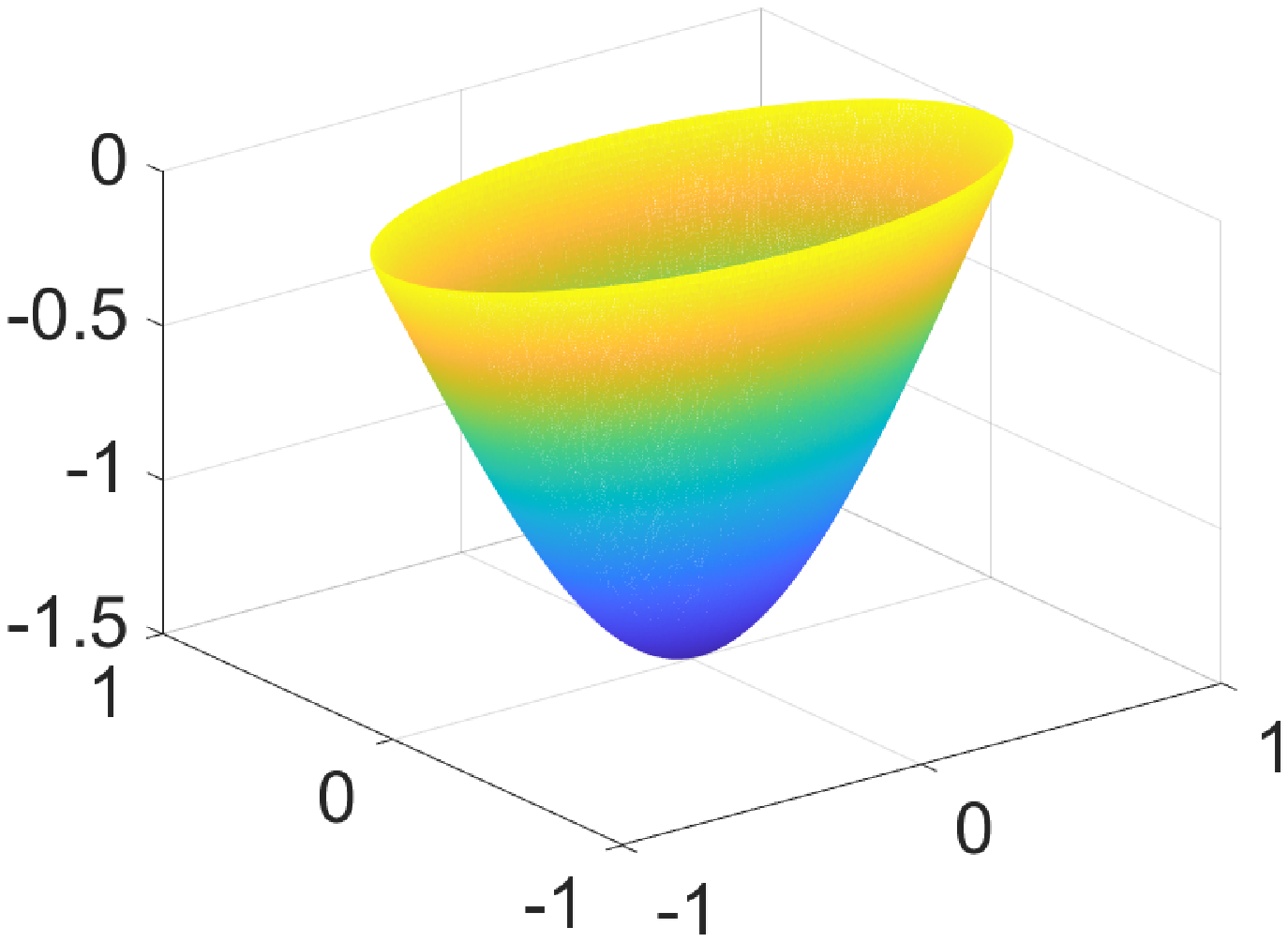}&
		\includegraphics[width=0.3\textwidth]{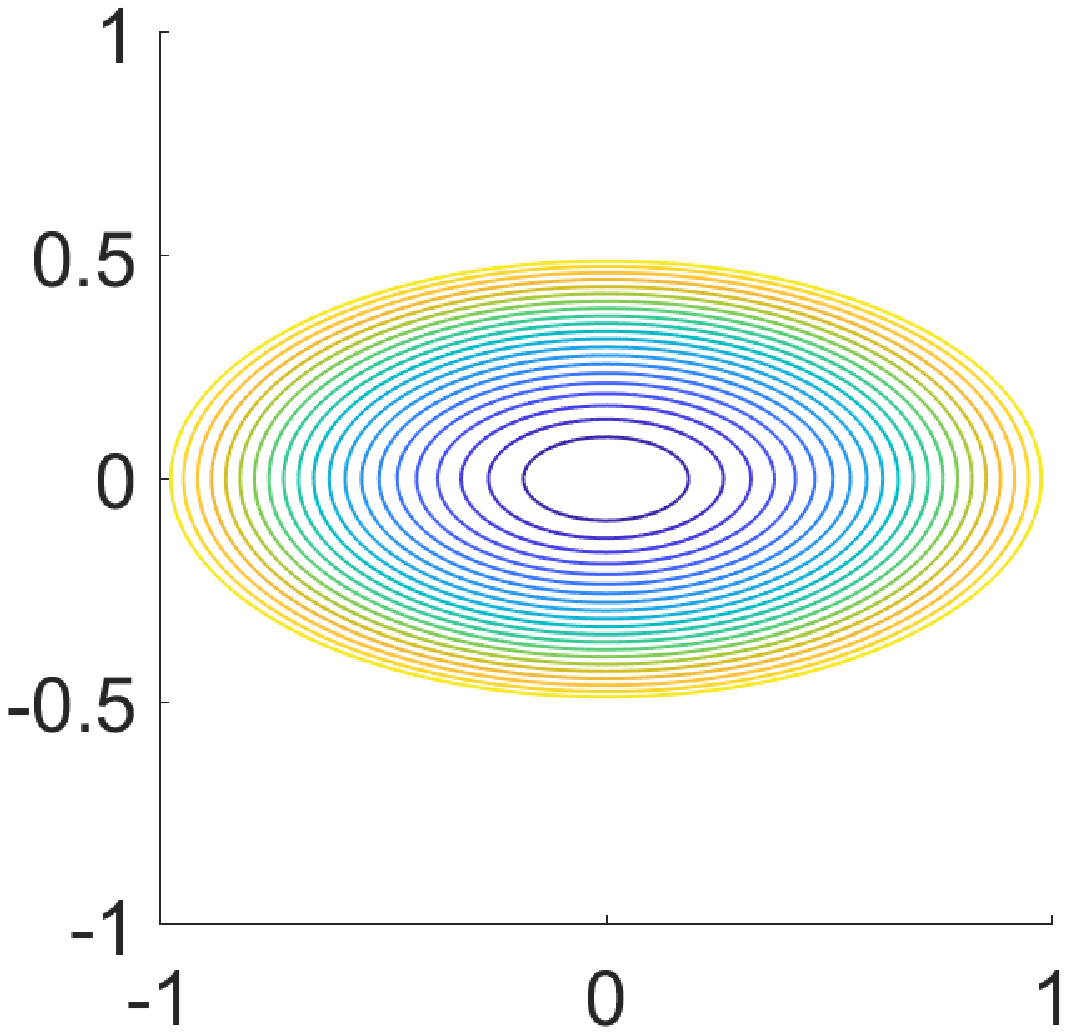}&
		\includegraphics[width=0.3\textwidth]{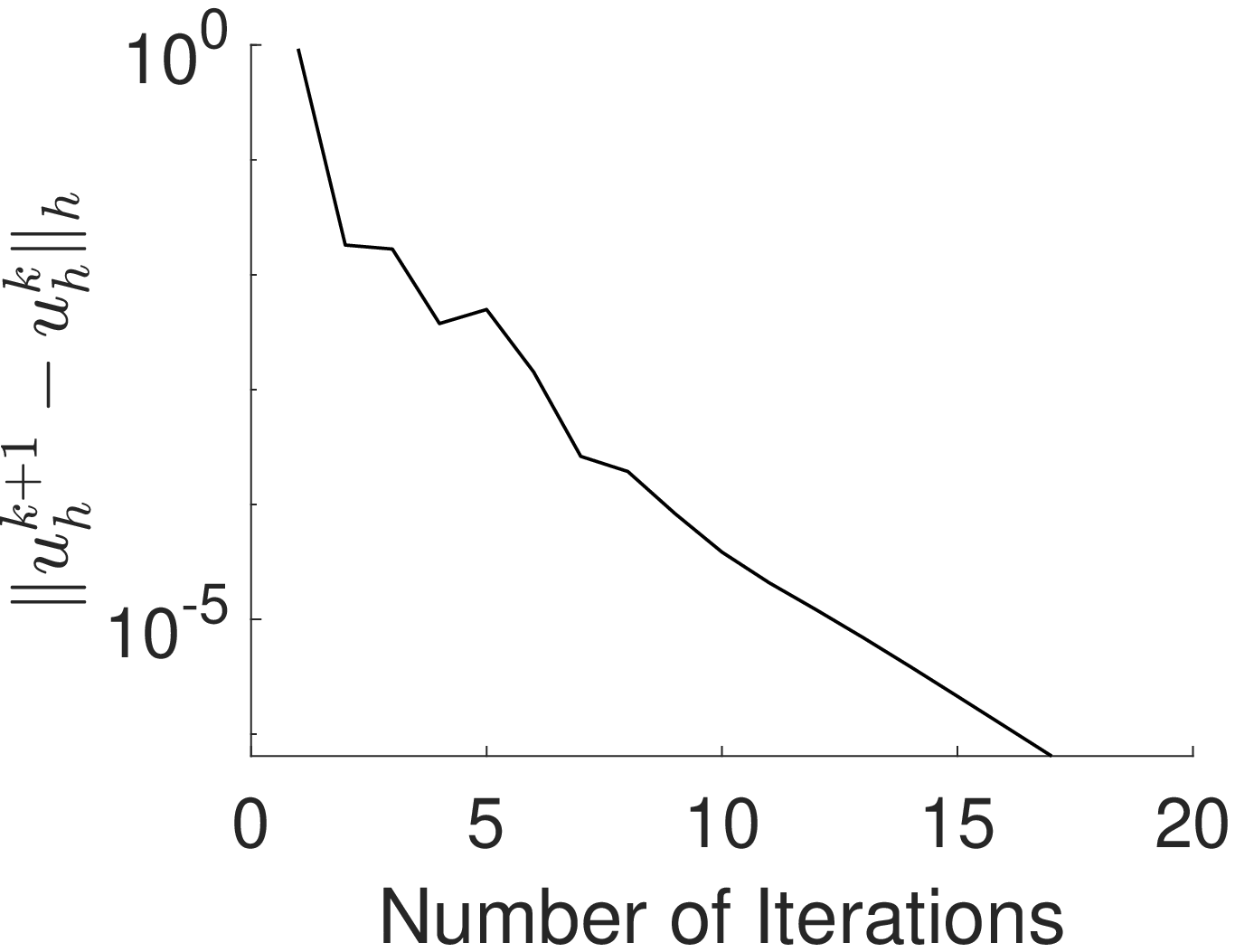} \\
		(d) & (e) & (f)\\
		\includegraphics[width=0.3\textwidth]{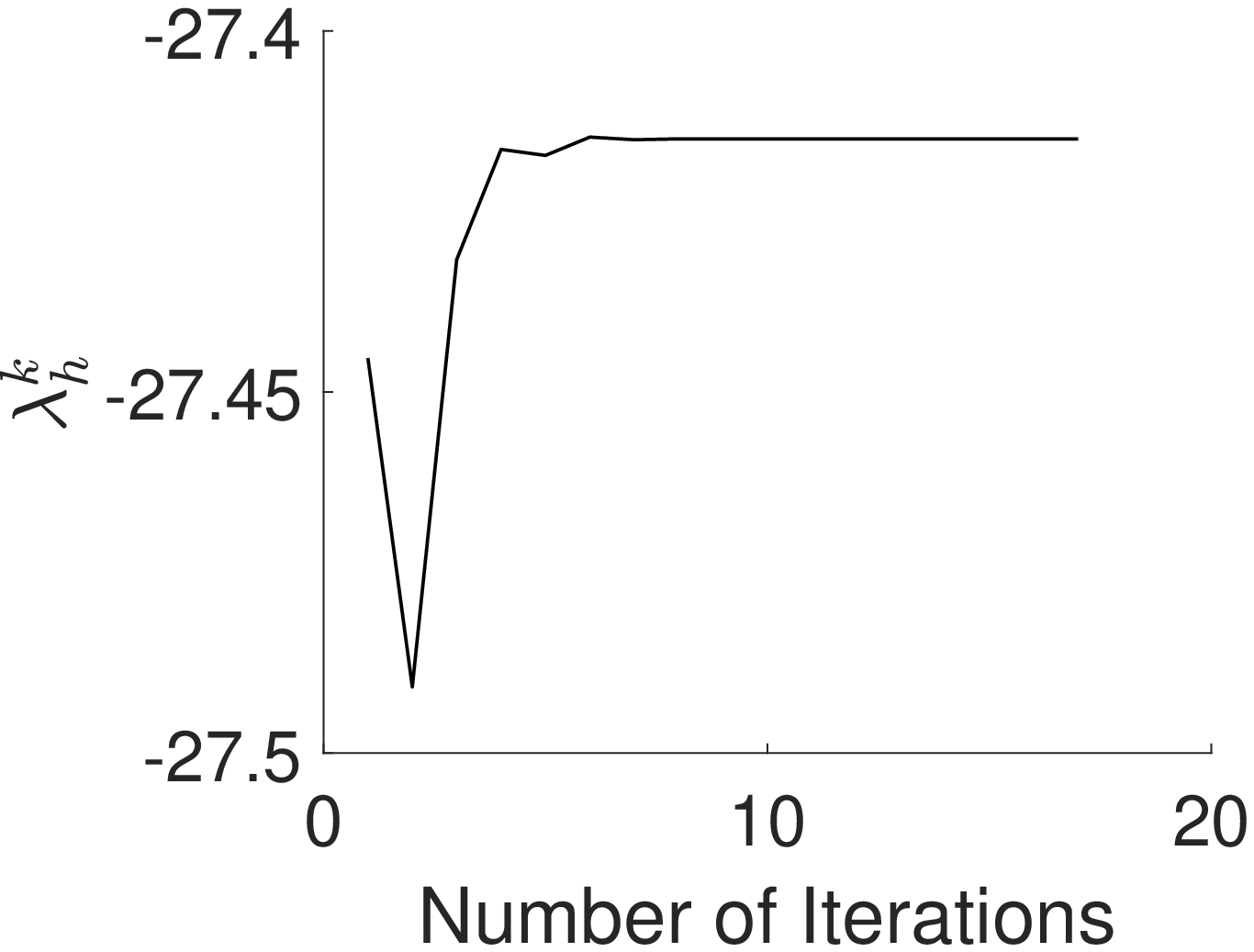}&
		\includegraphics[width=0.3\textwidth]{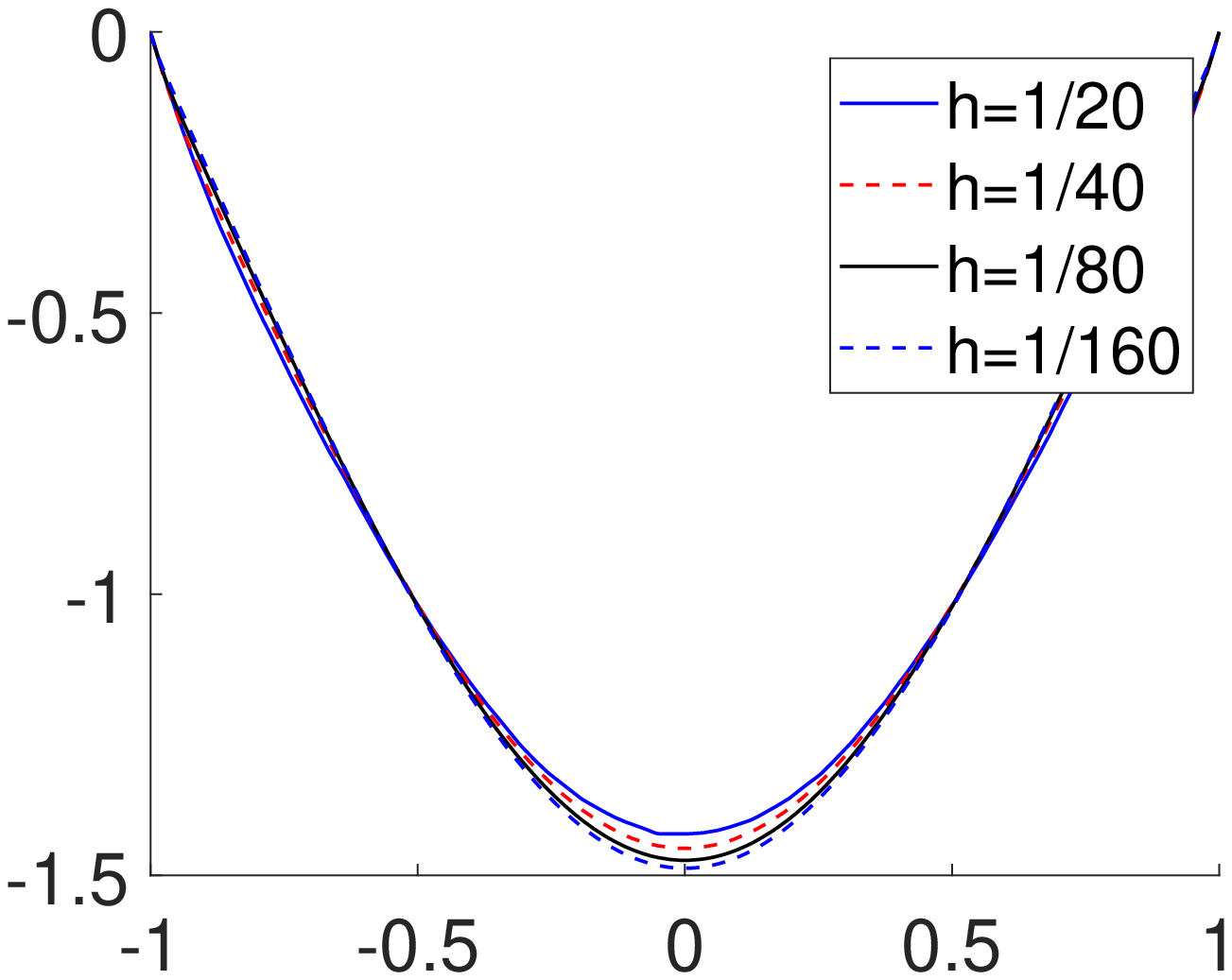}&
		\includegraphics[width=0.3\textwidth]{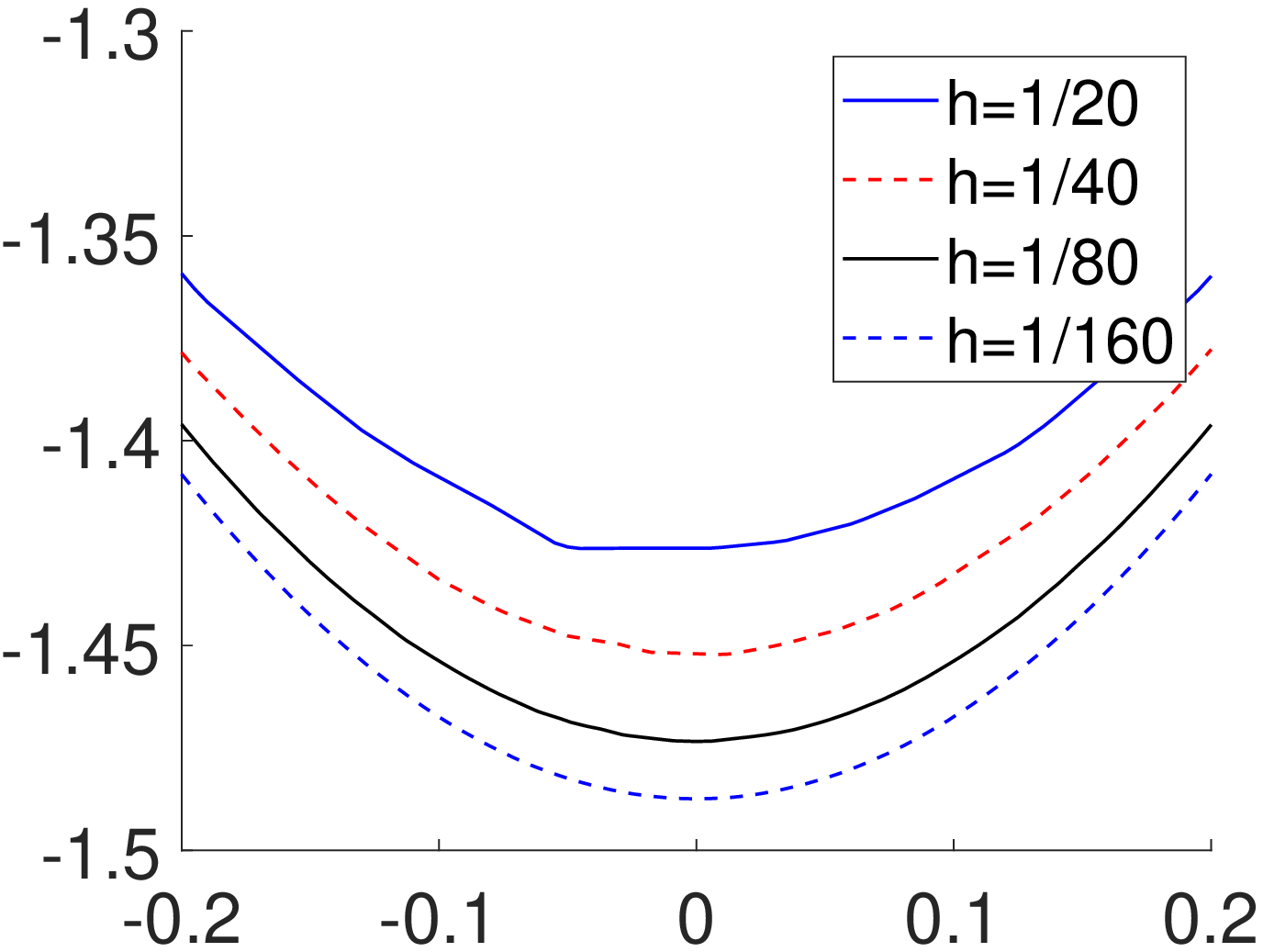}
	\end{tabular}
	\caption{The ellipse domain (\ref{eq.ellipse}). (a) The computed result with $h=1/80$. (b) The contour of (a). (c) The history of the error $\|u_h^{k+1}-u_h^k\|_h$. (d) The history of the computed eigenvalue $\lambda_h^k$. (e) Comparison of the cross sections along $x_2=0$ of the computed solution with various $h$. (f) Zoomed plot of the bottom region of (e).}
	\label{fig.ellipse}
\end{figure}

In the third example, we consider an ellipse domain defined by
\begin{align}
	\Omega=\left\{(x_1,x_2): x_1^2+2x_2^2<1\right\}.
	\label{eq.ellipse}
\end{align}
A triangulation of the domain with $h=1/20$ is visualized in Figure \ref{fig.domain}(c). In this set of experiments, we set stopping criterion $\xi=10^{-6}$ and time step $\tau=1/2$. The results with $h=1/80$ are shown in Figure \ref{fig.ellipse}(a)--(d). Similar to the results in the previous examples, the computed solution is smooth, and its contour consists of several ellipses with the same center, as shown in Figure \ref{fig.ellipse}(a) and Figure \ref{fig.ellipse}(b), respectively. In Figure \ref{fig.ellipse}(c), linear convergence is observed for the error $\|u_h^{k+1}-u_h^k\|_h$, and the convergence rate is about 0.34. The computed eigenvalue $\lambda_h^k$ attains its steady state with 6 iterations. With various $h$, we compare in Figure \ref{fig.ellipse}(e)--(f) the cross sections of the computed results along $x_2=0$. Convergence is observed as $h$ goes to 0.

With various $h$, the computational cost, the computed eigenvalue and minimal value of the computed solution are presented in Table \ref{tab.ellipse}. The eigenvalue $\lambda_h$ converges to $\lambda$ uniformly in the rate $\lambda_h\approx \lambda-ch$ with $\lambda\approx 29.5, c\approx 161$. In terms of the computational cost, all experiments used less than 20 iterations to satisfy the stopping criterion.


\begin{table}[t]
	\centering
	\begin{tabular}{|c|c|c|c|c|}
		\hline
		$h$& \# Iter.  & $\|u_h^{k+1}-u_h^k\|_{h}$ & $\lambda_h$& $\min u_h$\\
		\hline
		1/20 & 16 &6.80$\times 10^{-7}$ & 21.55 & -1.4277\\
		\hline
		1/40 & 16 & 9.68$\times 10^{-7}$ & 25.18 & -1.4525\\
		\hline
		1/80 & 17 & 6.44$\times 10^{-7}$ & 27.41 & -1.4734\\
		\hline
		1/160 & 17 & 7.00$\times 10^{-7}$ & 28.67 & -1.4875\\
		\hline
	\end{tabular}
	\caption{The ellipse domain (\ref{eq.ellipse}). Variations with $h$ of the number of iterations necessary to achieve convergence (2nd column), of the computed eigenvalue (4th column) and of the minimal value of $u_h$ over $\Omega$ (that is $u_h(\mathbf{0})$) (5th column).}
	\label{tab.ellipse}
\end{table}

\subsection{Example 4}
\begin{figure}
	\begin{tabular}{ccc}
		(a) & (b) & (c)\\
		\includegraphics[width=0.3\textwidth]{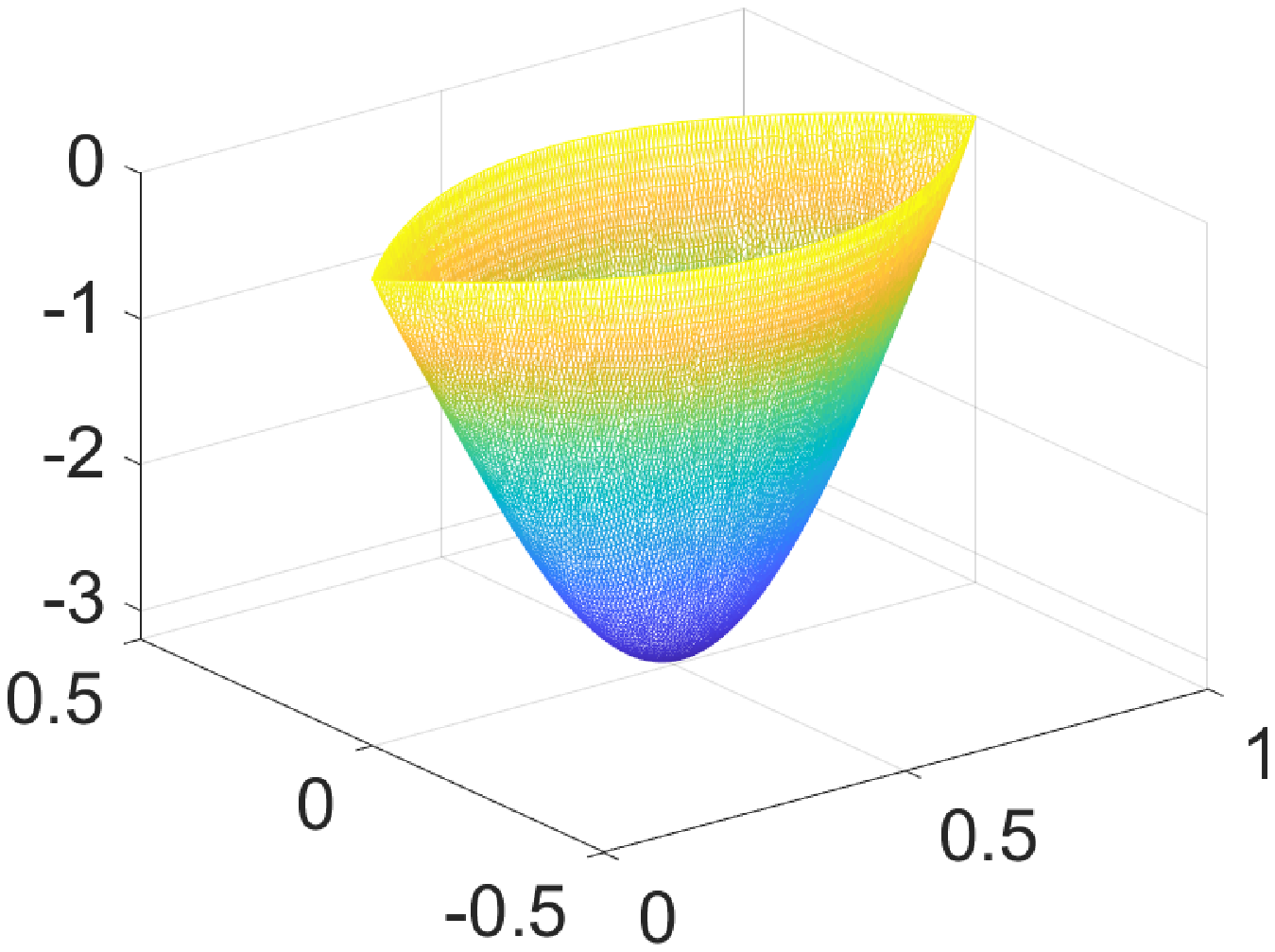}&
		\includegraphics[width=0.3\textwidth]{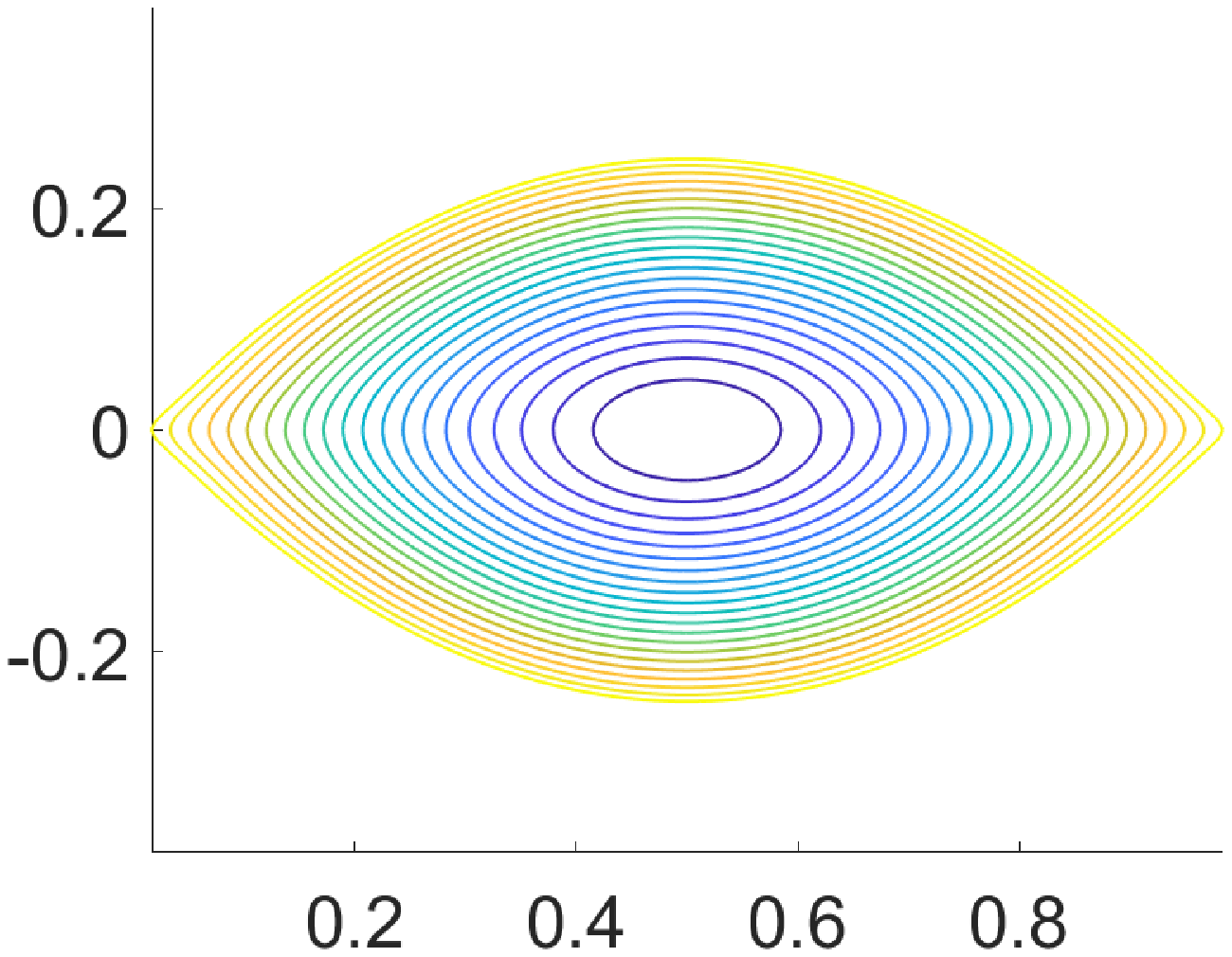}&
		\includegraphics[width=0.3\textwidth]{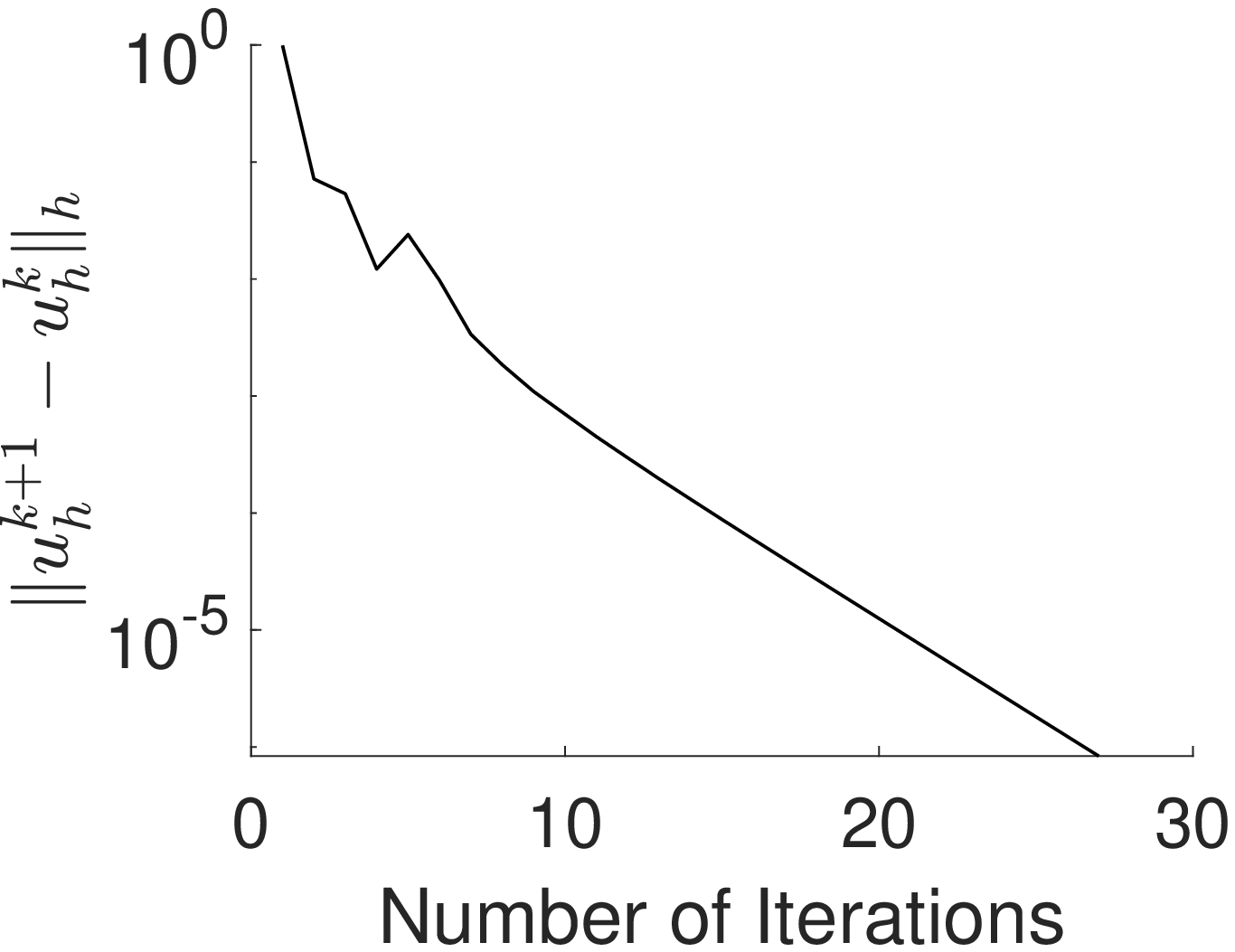} \\
		(d) & (e) & (f)\\
		\includegraphics[width=0.3\textwidth]{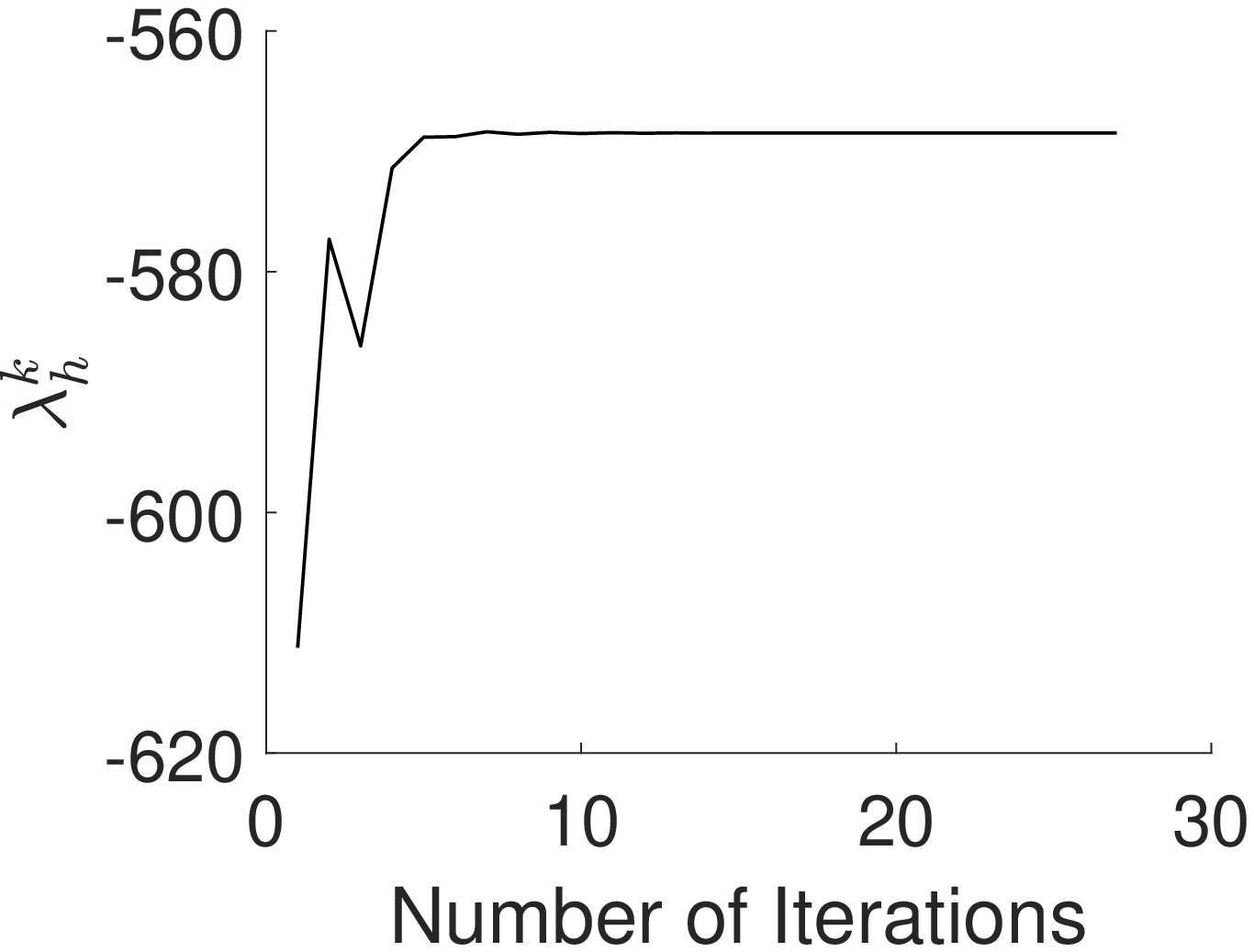}&
		\includegraphics[width=0.3\textwidth]{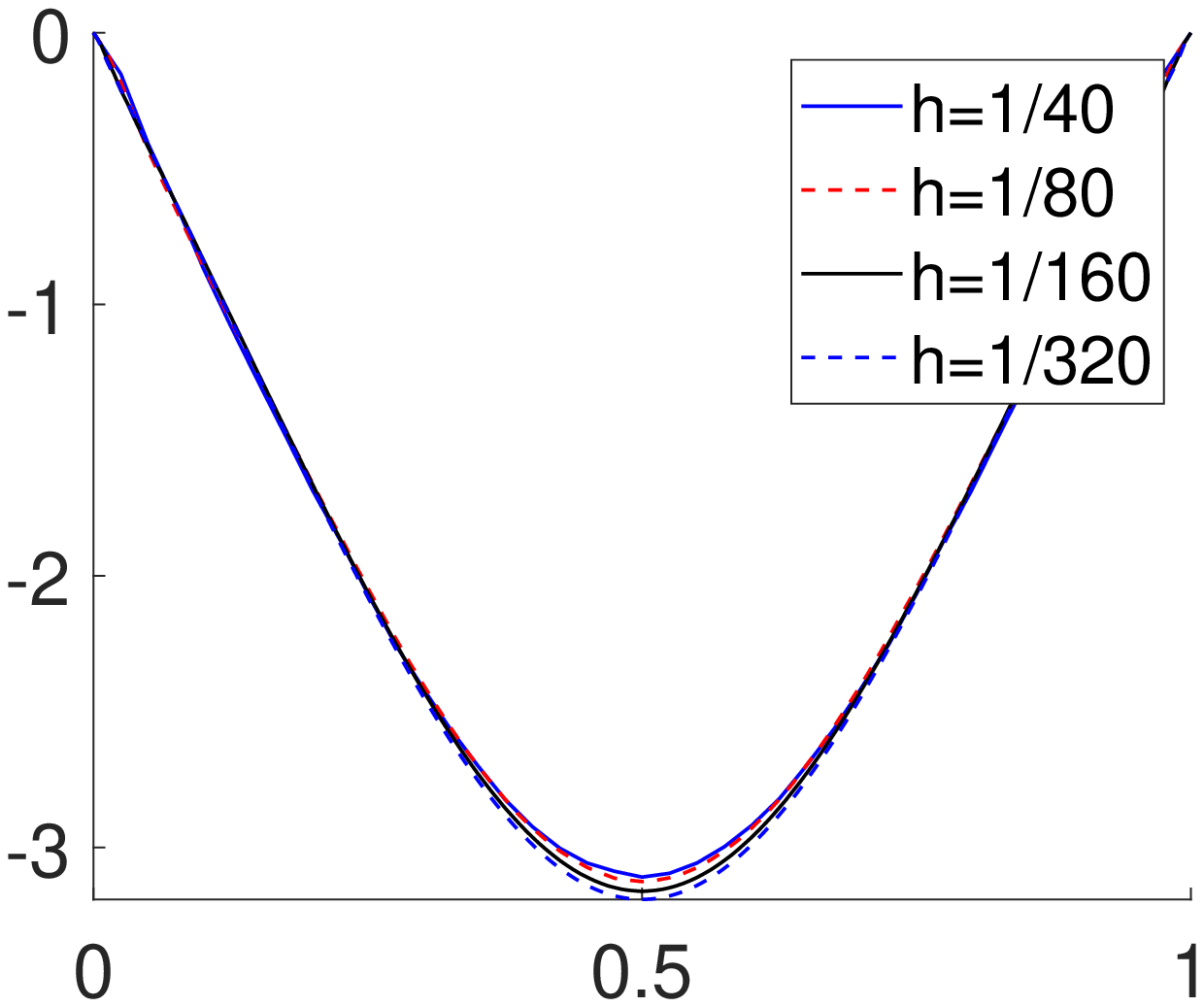}&
		\includegraphics[width=0.3\textwidth]{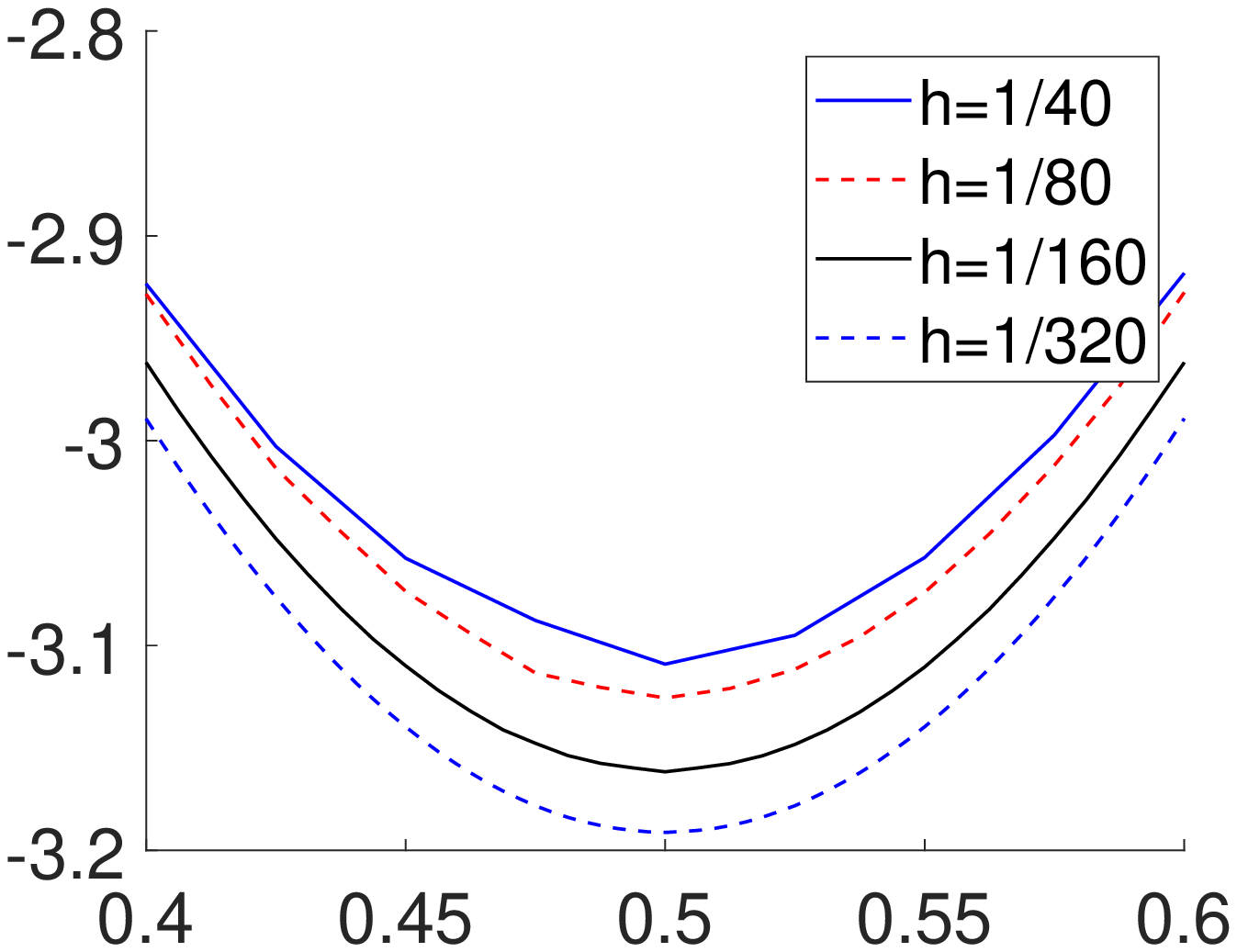}
	\end{tabular}
	\caption{The eye--shape domain (\ref{eq.eye}). (a) The computed result with $h=1/160$. (b) The contour of (a). (c)The history of the error $\|u_h^{k+1}-u_h^k\|_h$. (d) The history of the computed eigenvalue $\lambda_h^k$. (e) Comparison of the cross sections along $x_2=0$ of the computed solution with various $h$. (f) Zoomed plot of the bottom region of (e).}
	\label{fig.eye}
\end{figure}

We conclude this section by considering an open convex domain with a non-smooth boundary:
\begin{align}
	\Omega=\left\{ (x_1,x_2): -x_1(1-x_1)< x_2< x_1(1-x_1), \ 0<x_1<1\right\}.
	\label{eq.eye}
\end{align}
The domain described in the set (\ref{eq.eye}) has an eye shape, and its triangulation with $h=1/40$ is visualized in Figure \ref{fig.domain}(d). Since the domain is not smooth, in our experiments we use a smaller time step $\tau=1/8$ and larger regularization parameters $\varepsilon=4h^2$ and $c=4$. We set stopping criterion $\xi=10^{-6}$. The results with $h=1/160$  are shown in Figure \ref{fig.eye}(a)--(d). The computed solution is smooth, and its level curves have the same center, as shown in Figure \ref{fig.eye}(a) and Figure \ref{fig.eye}(b), respectively. In Figure \ref{fig.eye}(c), linear convergence is observed for the error $\|u_h^{k+1}-u_h^k\|_h$. The computed eigenvalue $\lambda_h^k$ attains its steady state with 7 iterations. With various $h$, we compare in Figure \ref{fig.eye}(e)--(f) the cross sections of the computed results along $x_2=0$. Convergence is observed as $h$ goes to 0.

With various $h$, the computational cost, the computed eigenvalue, and the minimal value of the computed solution are presented in Table \ref{tab.eye}. The eigenvalue $\lambda_h$ converges to $\lambda$ uniformly in the rate $\lambda_h\approx \lambda-ch$ with $\lambda\approx 618, c\approx 7792.3$. In terms of the computational cost, all experiments used no more than 30 iterations to satisfy the stopping criterion.

\begin{table}[t]
	\centering
	\begin{tabular}{|c|c|c|c|c|}
		\hline
		$h$& \# Iter.  & $\|u_h^{k+1}-u_h^k\|_{h}$ & $\lambda_h$& $\min u_h$\\
		\hline
		1/40 & 15 & 7.80$\times 10^{-7}$ & 425.51 & -3.1091\\
		\hline
		1/80 & 20 & 7.53$\times 10^{-7}$ & 516.57 & -3.1256\\
		\hline
		1/160 & 27 & 8.35$\times 10^{-7}$ & 568.47 & -3.1617\\
		\hline
		1/320 & 30 & 7.87$\times 10^{-7}$ & 597.39 & -3.1913\\
		\hline
	\end{tabular}
	\caption{The eye--shape domain (\ref{eq.eye}). Variations with $h$ of the number of iterations necessary to achieve convergence (2nd column), of the computed eigenvalue (4th column) and of the minimal value of $u_h$ over $\Omega$ (that is $u_h(\mathbf{0})$) (5th column).}
	\label{tab.eye}
\end{table}

\section{Conclusion}
We proposed an efficient operator--splitting method to solve the eigenvalue problem of the Monge--Amp\`{e}re equation. The backbone of our method relies on a convergent algorithm proposed in \cite{abedin2020inverse}. In each iteration, we solve a constrained optimization problem whose optimality condition is of the Monge--Amp\`{e}re type. We remove the constraint by including an indicator function and decouple the nonlinearity by introducing an auxiliary variable. The resulting problem is then converted to finding the steady state solution of an initial value problem which is time discretized by an operator--splitting method. The efficiency and effectiveness of the proposed method is demonstrated with several numerical experiments. In our experiments, we can choose a large constant time step. On smooth convex domains, our algorithm converges with a few iterations and is much faster than existing methods.
\label{sec.conclusion}

\section*{Acknowledgment}
Dr. Jun Kitagawa is acknowledged for bringing the Abedin-Kitagawa paper to the third author's attention in January 2021. Subsequently, Prof. Roland Glowinski and the authors initiated the project to develop an efficient algorithm to implement the Abedin-Kitagawa formulation.  

The preparation of this manuscript has been overshadowed by Roland's passing away in January 2022. Roland and the authors had intended to write jointly: most of the main ideas were worked out together and the authors have done their best to complete them. In sorrow, the authors dedicate this work to his memory. Roland's creativity, generosity, and friendship will be remembered. 

Hao Liu is partially supported by HKBU under grants 179356 and 162784.  S. Leung is supported by the Hong Kong RGC under grant 16302819. J. Qian is partially supported by NSF grants.
\bibliographystyle{abbrv}
\bibliography{reference}

\begin{thebibliography}{10}

\bibitem{abedin2020inverse}
F.~Abedin and J.~Kitagawa.
\newblock Inverse iteration for the {M}onge--{A}mp{\`e}re eigenvalue problem.
\newblock {\em Proceedings of the American Mathematical Society},
  148(11):4875--4886, 2020.

\bibitem{awanou2015standard}
G.~Awanou.
\newblock Standard finite elements for the numerical resolution of the elliptic
  {M}onge--{A}mp{\`e}re equation: classical solutions.
\newblock {\em IMA Journal of Numerical Analysis}, 35(3):1150--1166, 2015.

\bibitem{bakelman2012convex}
I.~J. Bakelman.
\newblock {\em Convex analysis and nonlinear geometric elliptic equations}.
\newblock Springer Science \& Business Media, 2012.

\bibitem{benamou2000computational}
J.-D. Benamou and Y.~Brenier.
\newblock A computational fluid mechanics solution to the {Monge-Kantorovich}
  mass transfer problem.
\newblock {\em Numerische Mathematik}, 84(3):375--393, 2000.

\bibitem{benamou2010two}
J.-D. Benamou, B.~D. Froese, and A.~M. Oberman.
\newblock Two numerical methods for the elliptic {M}onge--{A}mp{\`e}re
  equation.
\newblock {\em ESAIM: Mathematical Modelling and Numerical Analysis},
  44(4):737--758, 2010.

\bibitem{benamou2014numerical}
J.-D. Benamou, B.~D. Froese, and A.~M. Oberman.
\newblock Numerical solution of the optimal transportation problem using the
  {M}onge--{A}mp{\`e}re equation.
\newblock {\em Journal of Computational Physics}, 260:107--126, 2014.

\bibitem{bonito2016operator}
A.~Bonito, A.~Caboussat, and M.~Picasso.
\newblock Operator splitting algorithms for free surface flows: Application to
  extrusion processes.
\newblock In {\em Splitting Methods in Communication, Imaging, Science, and
  Engineering}, pages 677--729. Springer, 2016.

\bibitem{bukavc2013fluid}
M.~Buka{\v{c}}, S.~{\v{C}}ani{\'c}, R.~Glowinski, J.~Tamba{\v{c}}a, and
  A.~Quaini.
\newblock Fluid--structure interaction in blood flow capturing non-zero
  longitudinal structure displacement.
\newblock {\em Journal of Computational Physics}, 235:515--541, 2013.

\bibitem{caboussat2018least}
A.~Caboussat, R.~Glowinski, and D.~Gourzoulidis.
\newblock A least-squares/relaxation method for the numerical solution of the
  three-dimensional elliptic {M}onge--{A}mp{\`e}re equation.
\newblock {\em Journal of scientific computing}, 77(1):53--78, 2018.

\bibitem{caboussat2013least}
A.~Caboussat, R.~Glowinski, and D.~C. Sorensen.
\newblock A least-squares method for the numerical solution of the {D}irichlet
  problem for the elliptic {M}onge--{A}mp{\`e}re equation in dimension two.
\newblock {\em ESAIM: Control, Optimisation and Calculus of Variations},
  19(3):780--810, 2013.

\bibitem{caboussat2021second}
A.~Caboussat and D.~Gourzoulidis.
\newblock A second order time integration method for the approximation of a
  parabolic 2d {M}onge-{A}mp{\`e}re equation.
\newblock In {\em Numerical Mathematics and Advanced Applications ENUMATH
  2019}, pages 225--234. Springer, 2021.

\bibitem{caboussat2021adaptive}
A.~Caboussat, D.~Gourzoulidis, and M.~Picasso.
\newblock An adaptive method for the numerical solution of a 2d
  {M}onge-{A}mp{\`e}re equation.
\newblock In {\em Proceedings of the 10th International Conference on Adaptive
  Modeling and Simulation (ADMOS 2021)}, 2021.

\bibitem{dean2003numerical}
E.~J. Dean and R.~Glowinski.
\newblock Numerical solution of the two-dimensional elliptic
  {M}onge--{A}mp{\`e}re equation with dirichlet boundary conditions: an
  augmented {L}agrangian approach.
\newblock {\em Comptes rendus Math{\'e}matique}, 336(9):779--784, 2003.

\bibitem{dean2004numerical}
E.~J. Dean and R.~Glowinski.
\newblock Numerical solution of the two-dimensional elliptic
  {M}onge--{A}mp{\`e}re equation with dirichlet boundary conditions: a
  least-squares approach.
\newblock {\em Comptes rendus Math{\'e}matique}, 339(12):887--892, 2004.

\bibitem{deng2019new}
L.-J. Deng, R.~Glowinski, and X.-C. Tai.
\newblock A new operator splitting method for the {E}uler elastica model for
  image smoothing.
\newblock {\em SIAM Journal on Imaging Sciences}, 12(2):1190--1230, 2019.

\bibitem{engquist2016optimal}
B.~Engquist, B.~D. Froese, and Y.~Yang.
\newblock Optimal transport for seismic full waveform inversion.
\newblock {\em arXiv preprint arXiv:1602.01540}, 2016.

\bibitem{feng2013recent}
X.~Feng, R.~Glowinski, and M.~Neilan.
\newblock Recent developments in numerical methods for fully nonlinear second
  order partial differential equations.
\newblock {\em SIAM REVIEW}, 55(2):205--267, 2013.

\bibitem{feng2022narrow}
X.~Feng, T.~Lewis, and K.~Ward.
\newblock A narrow-stencil framework for convergent numerical approximations of
  fully nonlinear second order {PDE}s.
\newblock {\em arXiv preprint arXiv:2202.12782}, 2022.

\bibitem{feng2009mixed}
X.~Feng and M.~Neilan.
\newblock Mixed finite element methods for the fully nonlinear
  {M}onge--{A}mp{\`e}re equation based on the vanishing moment method.
\newblock {\em SIAM Journal on Numerical Analysis}, 47(2):1226--1250, 2009.

\bibitem{feng2009modified}
X.~Feng and M.~Neilan.
\newblock A modified characteristic finite element method for a fully nonlinear
  formulation of the semigeostrophic flow equations.
\newblock {\em SIAM Journal on Numerical Analysis}, 47(4):2952--2981, 2009.

\bibitem{feng2009vanishing}
X.~Feng and M.~Neilan.
\newblock Vanishing moment method and moment solutions for fully nonlinear
  second order partial differential equations.
\newblock {\em Journal of Scientific Computing}, 38(1):74--98, 2009.

\bibitem{froese2012numerical}
B.~D. Froese.
\newblock A numerical method for the elliptic {M}onge--{A}mp{\`e}re equation
  with transport boundary conditions.
\newblock {\em SIAM Journal on Scientific Computing}, 34(3):A1432--A1459, 2012.

\bibitem{froese2011convergent}
B.~D. Froese and A.~M. Oberman.
\newblock Convergent finite difference solvers for viscosity solutions of the
  elliptic {M}onge--{A}mp{\`e}re equation in dimensions two and higher.
\newblock {\em SIAM Journal on Numerical Analysis}, 49(4):1692--1714, 2011.

\bibitem{froese2011fast}
B.~D. Froese and A.~M. Oberman.
\newblock Fast finite difference solvers for singular solutions of the elliptic
  {M}onge--{A}mp{\`e}re equation.
\newblock {\em Journal of Computational Physics}, 230(3):818--834, 2011.

\bibitem{gilbarg1977elliptic}
D.~Gilbarg, N.~S. Trudinger, D.~Gilbarg, and N.~Trudinger.
\newblock {\em Elliptic Partial Differential Equations of Second Order}, volume
  224.
\newblock Springer, 1977.

\bibitem{glowinski2020numerical}
R.~Glowinski, S.~Leung, H.~Liu, and J.~Qian.
\newblock On the numerical solution of nonlinear eigenvalue problems for the
  {M}onge--{A}mp{\`e}re operator.
\newblock {\em ESAIM: Control, Optimisation and Calculus of Variations},
  26:118, 2020.

\bibitem{glowinski2015penalization}
R.~Glowinski, S.~Leung, and J.~Qian.
\newblock A penalization-regularization-operator splitting method for eikonal
  based traveltime tomography.
\newblock {\em SIAM Journal on Imaging Sciences}, 8(2):1263--1292, 2015.

\bibitem{glowinski2019finite}
R.~Glowinski, H.~Liu, S.~Leung, and J.~Qian.
\newblock A finite element/operator-splitting method for the numerical solution
  of the two dimensional elliptic {M}onge--{A}mp{\`e}re equation.
\newblock {\em Journal of Scientific Computing}, 79(1):1--47, 2019.

\bibitem{glowinski2017splitting}
R.~Glowinski, S.~J. Osher, and W.~Yin.
\newblock {\em Splitting Methods in Communication, Imaging, Science, and
  Engineering}.
\newblock Springer, 2017.

\bibitem{glowinski2016some}
R.~Glowinski, T.-W. Pan, and X.-C. Tai.
\newblock Some facts about operator-splitting and alternating direction
  methods.
\newblock In {\em Splitting Methods in Communication, Imaging, Science, and
  Engineering}, pages 19--94. Springer, 2016.

\bibitem{haker2004optimal}
S.~Haker, L.~Zhu, A.~Tannenbaum, and S.~Angenent.
\newblock Optimal mass transport for registration and warping.
\newblock {\em International Journal of computer vision}, 60(3):225--240, 2004.

\bibitem{he2020curvature}
Y.~He, S.~H. Kang, and H.~Liu.
\newblock Curvature regularized surface reconstruction from point clouds.
\newblock {\em SIAM Journal on Imaging Sciences}, 13(4):1834--1859, 2020.

\bibitem{kazdan1985prescribing}
J.~L. Kazdan.
\newblock {\em Prescribing the Curvature of a {R}iemannian Manifold},
  volume~57.
\newblock American Mathematical Soc., 1985.

\bibitem{le2017eigenvalue}
N.~Q. Le.
\newblock The eigenvalue problem for the {M}onge--{A}mp{\`e}re operator on
  general bounded convex domains.
\newblock {\em arXiv preprint arXiv:1701.05165}, 2017.

\bibitem{le2020convergence}
N.~Q. Le.
\newblock Convergence of an iterative scheme for the {M}onge--{A}mp{\`e}re
  eigenvalue problem with general initial data.
\newblock {\em arXiv preprint arXiv:2006.06564}, 2020.

\bibitem{lions1985two}
P.-L. Lions.
\newblock Two remarks on {M}onge--{A}mp{\`e}re equations.
\newblock {\em Annali di Matematica Pura ed Applicata}, 142(1):263--275, 1985.

\bibitem{liu2019finite}
H.~Liu, R.~Glowinski, S.~Leung, and J.~Qian.
\newblock A finite element/operator-splitting method for the numerical solution
  of the three dimensional {M}onge--{A}mp{\`e}re equation.
\newblock {\em Journal of Scientific Computing}, 81(3):2271--2302, 2019.

\bibitem{liu2021operator}
H.~Liu, X.-C. Tai, and R.~Glowinski.
\newblock An operator-splitting method for the {G}aussian curvature
  regularization model with applications in surface smoothing and imaging.
\newblock {\em arXiv preprint arXiv:2108.01914}, 2021.

\bibitem{liu2021color}
H.~Liu, X.-C. Tai, R.~Kimmel, and R.~Glowinski.
\newblock A color elastica model for vector-valued image regularization.
\newblock {\em SIAM Journal on Imaging Sciences}, 14(2):717--748, 2021.

\bibitem{liu2022elastica}
H.~Liu, X.-C. Tai, R.~Kimmel, and R.~Glowinski.
\newblock Elastica models for color image regularization.
\newblock {\em arXiv preprint arXiv:2203.09995}, 2022.

\bibitem{liu2022fast}
H.~Liu and D.~Wang.
\newblock Fast operator splitting methods for obstacle problems.
\newblock {\em arXiv preprint arXiv:2203.08380}, 2022.

\bibitem{roberts1995fully}
L.~A. Roberts, L.~A. Caffarelli, and X.~Cabr{\'e}.
\newblock {\em Fully Nonlinear Elliptic Equations}, volume~43.
\newblock American Mathematical Soc., 1995.

\bibitem{stojanovic2004optimal}
S.~Stojanovic.
\newblock Optimal momentum hedging via {M}onge--{A}mp{\`e}re {PDEs} and a new
  paradigm for pricing options.
\newblock {\em SIAM J. Control and Optimization}, 43:1151--1173, 2004.

\bibitem{tso1990real}
K.~Tso.
\newblock On a real {M}onge--{A}mp{\`e}re functional.
\newblock {\em Inventiones Mathematicae}, 101(1):425--448, 1990.

\end{thebibliography}
\end{document}